\documentclass[a4j,11pt]{article}

\usepackage{graphicx}
\usepackage{framed}
\usepackage{latexsym}
\usepackage{amsthm}
\usepackage{bm}
\usepackage{natbib}
\newtheorem{theorem}{Theorem}
\newtheorem{proposition}{Proposition}
\newtheorem{corollary}{Corollary}
\newtheorem{lemma}{Lemma}

\def\vec#1{\mbox{\boldmath $#1$}}
\usepackage{amsmath,amssymb}

\usepackage{amssymb}
\usepackage{latexsym}
\def\vec#1{\mbox{\boldmath $#1$}}
\usepackage{amsmath}

\DeclareMathOperator{\diag}{diag}

\DeclareMathOperator*{\argmin}{argmin}

\begin{document}

\title{Asymptotics for penalized splines in generalized additive models}

\author{
{\sc Takuma Yoshida}$^{1}$\ \ {\sc and}\ \ {\sc Kanta Naito}$^{1}$\thanks{E-mail:\ {\tt naito@riko.shimane-u.ac.jp}}\\
$^{1}${\it Graduate School of Science and Engineering} \\
{\it Shimane University, Matsue 690-8504, Japan}
}

\date{\empty}

\maketitle

%
\begin{abstract}                                             
This paper discusses asymptotic theory for penalized spline estimators in generalized additive models. 
The purpose of this paper is to establish the asymptotic bias and variance as well as the asymptotic normality of the penalized spline estimators proposed by Marx and Eilers (1998). 
Furthermore, the asymptotics for the penalized quasi likelihood fit in mixed models are also discussed. 
\end{abstract}


\section{Introduction}

\quad 
The generalized additive model(GAM) is a typical regression model, in which the relationship between the one-dimensional response $Y$ and the multidimensional explanatory $\vec{x}=(x_1,\cdots,x_D)$ is modeled by a link function $g(\cdot)$, as follows:
$$
g(E[Y|\vec{X}=x])=\eta(\vec{x})=\eta_1(x_1)+\cdots+\eta_D(x_D),
$$
where each $\eta_j(j=1,\cdots,D)$ is a univariate regression function. 
If $Y$ has a Gaussian distribution, then $g$ is the identity function and, hence, $\eta(\vec{x})=E[Y|\vec{X}=\vec{x}]$. 
Additionally, the GAM can specify a distribution such as a Bernoulli, Poisson or Gamma distribution. 
In GAMs, the purpose is often to estimate $\eta$.
The parametric and the nonparametric estimation techniques of $\eta$ have been established by several authors (see Hastie et al. (1990) and Wood (2006)). 
In this paper, $\eta$ is estimated via the penalized spline method. 
Penalized splines were introduced by O'Sullivan (1986) and Eilers and Marx (1996), and are recognized as an efficient technique for GAMs.
Applications and theories of penalized splines in GAMs have been widely discussed, including by Aerts et al. (2002) and Ruppert et al. (2003).

To construct the estimator of $\eta_j$'s, a repetition update method, the so-called backfitting algorithm, is often used. 
However, when the response has a non-Gaussian distribution, such as a Bernoulli or Poisson distribution, 
the overall estimation procedure becomes complicated and its computation time grows, because the estimators are obtained by using a blend of backfitting and the Fisher-scoring algorithm. 
On the other hand, Marx and Eilers (1998) proposed a new penalized spline estimator without backfitting algorithms, which is denoted as the ridge corrected penalized spline estimator (RCPS). 
We will briefly describe the RCPS as we will focus on it later. 
The penalized spline estimator is obtained based on maximization of the penalized log-likelihood. 
However, it appears difficult to obtain the maximizer of the penalized log-likelihood $\ell$ as the Hessian of $\ell$ is not invertible.  
The RCPS is constructed based on the maximization of $\ell_\gamma$, which is $\ell$ plus an additional ridge penalty. 
Since the Hessian of $\ell_\gamma$ is invertible, the maximizer of $\ell_\gamma$ can be obtained via the Fisher-scoring algorithm. 
Thus, it is easy to construct the RCPS.

In univariate models($D=1$),  Hall and Opsomer (2005), Claeskens et al. (2009),
Kauermann et al. (2009) and Wang et al. (2011) researched the asymptotic properties of penalized spline estimators. 
Recently, Yoshida and Naito (2012) worked on the asymptotic distribution of penalized splines in an additive model.
In contrast, Horowitz and Mammen (2004), Linton (2000) and Yu et al. (2008) studied the asymptotics for the kernel estimator in a GAM. 
However, the asymptotic results for penalized spline estimators in GAMs have not yet been sufficiently developed like they have been for their applications. 

In this paper, the asymptotics for penalized splines in a GAM are discussed. 
Our main purpose is to establish the asymptotic normality of the RCPS. 
Kauermann et al. (2009) showed the asymptotic normality of the penalized spline estimator in generalized linear models(GLM). 
Hence, the results in this paper generalize the results of Kauermann et al. (2009).
Furthermore, penalized spline smoothing can be linked to mixed models(see Lin and Zhang (1999) and Ruppert et al. (2003)). 
In generalized additive mixed models (GAMM), the penalized quasi likelihood method (PQL) is an efficient method for obtaining the estimator and predictor. 
We also show the asymptotic normality of the PQL fit.

This paper is organized as follows. 
In Section 2, the GAM is defined and the RCPS is constructed by the Fisher-scoring algorithm. 
Section 3 shows the asymptotic normality of the RCPS. 
Section 4 states the asymptotics for the PQL fits in a GAMM. 
In Section 5, the applications of the approximate confidence interval are addressed. 
Section 6 provides a numerical study to validate the asymptotic normality of the RCPS. 
Related discussions and issues for future research are addressed in Section 7, and proofs for theoretical results are given in the Appendix.


\section{Penalized spline estimator in a GAM}

\subsection{Generalized additive spline model}

For the dataset $\{ (y_i,\vec{x}_{i})|i=1,\cdots,n\}$, consider an exponential family of the generalized additive model with a canonical link function 
\begin{eqnarray}
f(y_i|\vec{x}_i,\eta)=\exp \left( \frac{y_i\eta(\vec{x}_i)-c(\eta(\vec{x}_i))}{\phi } + h(y_i,\phi) \right),\quad i=1,\cdots,n, \label{model1}
\end{eqnarray}
where $\vec{x}_i=(x_{i1},\cdots,x_{iD})$ is a $D$-variate explanatory variable, $\eta(\cdot)$ is an unknown natural parameter which has the additive formation 
$$
\eta(\vec{x})=\eta_1(x_1)+\cdots+\eta_D(x_D),
$$ 
where $\eta_j$'s is an unknown univariate function, $\phi$ is a dispersion parameter, and $c$ and $h$ are known functions. 
The canonical link function indicates $g^{-1}=c^\prime$, which leads to 
$E[Y_i|\vec{X}_i=\vec{x}_{i}]=g^{-1}(\eta(\vec{x}_i))=c^\prime(\eta(\vec{x}_i))$ and $V[Y_i |\vec{X}_i=\vec{x}_{i}]= \phi c^{\prime\prime}(\eta(\vec{x}_i ))$, where $c^\prime$ and $c^{\prime\prime}$ are the first and second derivatives of $c$. 
More general settings concerned with the link function were clarified by McCullagh and Nelder (1989).
In this paper, for the natural parameter $\eta(\vec{x})$, we assume that $E[\eta_j(X_j)]=0 (j=1,\cdots,D)$ to ensure the identifiability of $\eta_j$. 
For simplicity, we hereafter ignore the role of the dispersion parameters in (\ref{model1}) and set $\phi\equiv 1$, thus denoting $h(y,\phi)=h(y)$. 

Our purpose is to estimate $\eta_j$ using nonparametric spline methods.
We now prepare the $B$-spline model
$$
s(x)=\sum_{k=-p+1}^{K_{n}} B_k^{[p]}(x)b_{k,j},\ \ j=1,\cdots,D
$$ 
as an approximation to $\eta_j(x)$, 
where $B_k^{[p]}(x)(k=-p+1,\cdots,K_{n})$ are the $p$th $B$-spline functions defined recursively as 
\begin{eqnarray*}
B_k^{[0]}(x)&=&
\left\{
\begin{array}{cc}
1,& \kappa_{k-1}<x\leq \kappa_k,\\
0,& {\rm otherwise},
\end{array}
\right. \\
B_k^{[p]}(x)&=&\frac{x-\kappa_{k-1}}{\kappa_{k+p-1}-\kappa_{k-1}}B_k^{[p-1]}(x)+\frac{\kappa_{k+p}-x}{\kappa_{k+p}-\kappa_{k}}B_{k+1}^{[p-1]}(x),
\end{eqnarray*}
where $\kappa_k (k=-p+1,\cdots,K_n+p)$ are knots and $b_{k,j}$'s is an unknown parameter. 
Some fundamental properties of $B$-splines were detailed by de Boor (2001). 
We will write $B_k^{[p]}(x)=B_k(x)$ unless we specify the degree of $B$-splines, because we will mainly focus on the $p$th $B$-spline from now on. 
The suggested density of $Y_i$ is defined as 
\begin{eqnarray*}
f(y_i|\vec{x}_i,\vec{b})=\exp \left( y_i Z(\vec{x}_i)^T\vec{b}-c(Z(\vec{x}_i)^T\vec{b}) + h(y_i) \right),
\end{eqnarray*}
where $Z(\vec{x})=(\vec{B}(x_1)^T\ \cdots\ \vec{B}(x_D)^T)^T$, $\vec{B}(x)=(B_{-p+1}(x)\ \cdots\ B_{K_n}(x))^T$, $\vec{b}=(\vec{b}_1^T\ \cdots\ \vec{b}_D^T)^T$ and $\vec{b}_j=(b_{-p+1,j}\ \cdots\ b_{K_n,j})^T$. 
Using the estimator $\hat{\vec{b}}_j=(\hat{b}_{-p+1,j}\ \cdots\ \hat{b}_{K_n,j})^T$ of $\vec{b}_j$, the estimator of $\eta_j(x_j)$ can be defined as 
$$
\hat{\eta}_j(x_j)=\sum_{k=-p+1}^{K_n} B_{k}(x_j)\hat{b}_{k,j},\ \ j=1,\cdots,D.
$$

\subsection{The penalized spline estimator}

To estimate $\vec{b}$, we prepare the log-likelihood 
\begin{eqnarray*}
\frac{1}{n}\sum_{i=1}^n \log f(y_i|\vec{x}_i,\vec{b})=\frac{1}{n}\{\vec{y}^T (Z\vec{b})-\vec{1}^T c(Z\vec{b})\}+\frac{1}{n}\vec{1}^T h(\vec{y}),
\end{eqnarray*} 
where $\vec{y}=(y_1\ \cdots\ y_n)^T$, $c(Z\vec{b})=(c(Z(\vec{x}_1)^T\vec{b})\ \cdots\ c(Z(\vec{x}_n)^T\vec{b}))^T$, $Z_k=(B_{-p+j}(x_{ik}))_{ij}$, $Z=[Z_1\ \cdots\ Z_D]$, and $h(\vec{y})=(h(y_1),\cdots,h(y_n))^T$. 
It is known that the spline estimator obtained by maximization of the log-likelihood tends to display 'wiggle behavior'. 
Hence, we consider using the penalized spline estimator to obtain a smooth curve. 
Define the penalized log-likelihood as follows
\begin{eqnarray}
\ell(\vec{b},\lambda_n)
&=&
\frac{1}{n}\sum_{i=1}^n \log f(y_i|\vec{x}_i,\vec{b})-\sum_{j=1}^D \frac{\lambda_{jn}}{2n}\vec{b}_j^T \Delta_m^\prime \Delta_m \vec{b}_j\nonumber\\
&=&
\frac{1}{n}\{\vec{y}^T (Z\vec{b})-\vec{1}^T c(Z\vec{b})\}+\frac{1}{n}\vec{1}^T h(\vec{y})-\frac{1}{2n} \vec{b}^T Q_m(\lambda_n) \vec{b}, \label{pen1}
\end{eqnarray}
where $\lambda_{jn}$ is the smoothing parameter $(j=1,\cdots,D)$, the $(K_n+p-m)\times (K_n+p)$th matrix $\Delta_m$ is the $m$th difference matrix, which is given by Marx and Eilers (1998) and $Q_m(\lambda_n)={\rm diag}[\lambda_{1n}\Delta_m^\prime \Delta_m\ \cdots\ \lambda_{Dn}\Delta_m^\prime \Delta_m]$. 
In general, the maximizer of (\ref{pen1}) is obtained by the Fisher-scoring algorithm.  
As in the typical problem of spline methods in a GAM, however, the Hessian of $\ell(\vec{b},\lambda_n)$ is not invertible and so the Fisher-scoring method is not usable directly. 
To overcome this singularity problem, we can use backfitting algorithms (see Hastie and Tibshirani (1990)). 
However, the overall algorithm becomes complicated and the computation grows (see Section 1). 
These problems were discussed by Marx and Eilers (1998).
We will next review the ridge corrected penalized spline estimator.

\subsection{The ridge corrected penalized spline estimator}

Marx and Eilers (1998) proposed a nice estimation method for $\vec{b}$ without using a backfitting algorithm for penalized splines in the GAM context.
They defined the ridge corrected penalized log-likelihood as
\begin{eqnarray}
\ell(\vec{b},\lambda_n,\gamma_n)
=
\ell(\vec{b},\lambda_n)-\frac{\gamma_n}{2n}\vec{b}^T \vec{b}, \label{pen}
\end{eqnarray}
where $\gamma_n>0$. 
Let $\hat{\vec{b}}=(\hat{\vec{b}}_1^T\ \cdots\ \hat{\vec{b}}_D^T)^T$ be the maximizer of (\ref{pen}), which can be obtained directly via the Fisher-scoring method since the Hessian of $\ell(\vec{b},\lambda_n,\gamma_n)$ is invertible. 
The gradient $G(\vec{b},\lambda_n,\gamma_n)$ and Hessian $H(\vec{b},\lambda_n,\gamma_n)$ of $\ell(\vec{b},\lambda_n,\gamma_n)$ are obtained with
\begin{eqnarray*}
G(\vec{b},\lambda_n,\gamma_n)&=&\frac{\partial \ell(\vec{b},\lambda_n,\gamma_n)}{\partial \vec{b}}
=
\frac{1}{n}\{Z^T\vec{y}-Z^T c^\prime (Z\vec{b})\}-\frac{1}{n}Q_m(\lambda_n)\vec{b}-\frac{\gamma_n}{n}\vec{b},\\
H(\vec{b},\lambda_n,\gamma_n)&=&\frac{\partial^2 \ell(\vec{b},\lambda_n,\gamma_n)}{\partial \vec{b}\partial \vec{b}^T}
=
-\frac{1}{n}Z^T \diag[c^{\prime\prime} (Z\vec{b})]Z-\frac{1}{n}Q_m(\lambda_n)-\frac{\gamma_n}{n} I,
\end{eqnarray*} 
where $c^{\prime}(Z\vec{b})$ and $c^{\prime\prime}(Z\vec{b})$ are $n$-vectors  defined in the same manner as $c(Z\vec{b})$.
The $k$-step iterated estimator $\vec{b}^{(k)}$ of $\vec{b}$ can be written as
\begin{eqnarray*}
{\small
\vec{b}^{(k)}=
(Z^T W^{(k-1)} Z+Q_m(\lambda_n)+\gamma_n I)^{-1}
Z^T W^{(k-1)}\left\{ Z\vec{b}^{(k-1)}+(W^{(k-1)})^{-1}\left\{\vec{y}-c^\prime (Z\vec{b}^{(k-1)})\right\}\right\},
}
\end{eqnarray*}
where $W^{(k-1)}={\rm diag}[c^{\prime\prime}(Z\vec{b}^{(k-1)})]$.
As $k\rightarrow \infty$, $\vec{b}^{(k)}$ converges to $\hat{\vec{b}}$ if the initial $\vec{b}^{(0)}$ is appropriately chosen. 
The RCPS of $\eta_j(x_j)$ can be obtained as $\hat{\eta}_j(x_j)=\vec{B}(x_j)^T \hat{\vec{b}}_j$.
In the next section, we discuss the asymptotic properties of $[\hat{\eta}_1(x_1)\ \cdots\ \hat{\eta}_D(x_D)]^T$.

\section{Asymptotic theory}

Here, we list some assumptions regarding the asymptotics of the penalized spline estimator.
\\

{\noindent \bf Assumptions}
\begin{enumerate}
\item[1.] The explanatory $\vec{X}=(X_1,\cdots,X_D)$ is distributed as $P(\vec{x})$ on $[0,1]^D$, where $[0,1]^D$ is the $D$-variate unit cube. 
\item[2.] For $j=1,\cdots,D$, $\eta_j\in C^{p+1}$ and $c\in C^3$. 
\item[3.] The knots for the $B$-spline basis are equidistantly located with $\kappa_k=k/K_n(k=-p+1,\cdots,K_n+p)$ and the number of knots satisfies $K_n=o(n^{1/2})$. 
\item[4.] For the non-singularity of $H(\vec{b},\lambda_n,\gamma_n)$, $K_n$ is chosen such that $D(K_n+p)<n$.
\item[5.] The smoothing parameters $\lambda_{jn} (j=1,\cdots,D)$ are positive sequences such that $\lambda_{jn}^{-1}$ is larger than the maximum eigenvalue of $(Z_j^T Z_j)^{-1/2}\Delta_m^\prime \Delta_m(Z_j^T Z_j)^{-1/2}$. 
\item[6.] Lastly, $\gamma_n=o(\lambda_nK_n^{-m})$, where $\lambda_n=\max_{j}\{\lambda_{jn}\}$.
\end{enumerate}

For a random variable $U_n$, $E[U_n|\vec{X}_n]$ and $V[U_n|\vec{X}_n]$ denote the conditional expectation and variance of $U_n$ given $(\vec{X}_1,\cdots,\vec{X}_n)=(\vec{x}_1,\cdots,\vec{x}_n)$, respectively. 
Define the $(K_n+p)$th square matrix $G_{k}=(G_{k,ij})_{ij}$, where the $(i,j)$-th component is 
\begin{eqnarray*}
G_{k,ij}=\int_{[0,1]^D} c^{\prime\prime}(\eta(\vec{x}))B_{-p+i}(x_k)B_{-p+j}(x_k)dP(\vec{x}).
\end{eqnarray*} 
Using this, we get $\Gamma_j(\lambda_{jn})=(G_{j}+(\lambda_{jn}/n)\Delta_m^\prime \Delta_m)$.
Let  
\begin{eqnarray}
\vec{b}_0=(\vec{b}_{10}^T \cdots\ \vec{b}_{D0}^T)^T=\underset{\vec{b}}{\argmin}\left\{\frac{1}{n}\sum_{i=1}^nE\left[\left.\log \frac{f(Y_i|\vec{x}_i,\eta)}{f(Y_i|\vec{x}_i,\vec{b})}\right|\vec{X}_n\right]\right\} \label{b0}
\end{eqnarray}
and let $\eta_{j0}(x_j)=\vec{B}(x_j)^T \vec{b}_{j0} (j=1,\cdots,D)$.

The asymptotic bias of $\hat{\eta}_j(x_j)$ can be written as 
$$
E[\hat{\eta}_j(x_j)|\vec{X}_n]-\eta_j(x_j)=E[\hat{\eta}_j(x_j)|\vec{X}_n]-\eta_{j0}(x_j)+\eta_{j0}(x_j)-\eta_j(x_j).
$$
In the following Proposition \ref{app}, the difference $\eta_{j0}(x_j)-\eta_j(x_j)$ is asymptotically evaluated. 
The asymptotics for $\hat{\eta}_j(x_j)-\eta_{j0}(x_j)=\vec{B}(x_j)^T (\hat{\vec{b}}_j-\vec{b}_{j0})$ can be shown in the following Theorem \ref{mv} by using the Taylor expansion of $G(\hat{\vec{b}},\lambda_n,\gamma_n)$ around $\vec{b}_0$ (see Lemma \ref{ex} in the Appendix), the properties of a partitioned matrix of $H(\vec{b},\lambda_n,\gamma_n)$ and its asymptotic results. 

\begin{proposition}\label{app}
Under the Assumptions, for $j=1,\cdots,D,$
\begin{eqnarray*}
\eta_{j0}(x_j)-\eta_j(x_j)= b_{j,a}(x_j)+o(K_n^{-(p+1)}),
\end{eqnarray*}
where 
$$
b_{j,a}(x)=-\frac{\eta^{(p+1)}_j(x)}{K_n^{p+1} (p+1)!}\sum_{k=1}^{K_n}I(\kappa_{k-1}\leq x<\kappa_k){\rm Br}_{p+1}\left(\frac{x-\kappa_{k-1}}{K_n^{-1}}\right),
$$
$I(a<x<b)$ is the indicator function of an interval $(a,b)$ and ${\rm Br}_p(x)$ is the $p$th Bernoulli polynomial.
\end{proposition}

\begin{theorem}\label{mv} 
Under the Assumptions, for $j=1,\cdots,D,$
\begin{eqnarray*}
E[\hat{\eta}_j(x_j)|\vec{X}_n]-\eta_j(x_j)
&=&b_{j,a}(x_j)+b_{j,\lambda}(x_j)+o_P(K_n^{-(p+1)})+o_P(\lambda_{jn}K_n^{1-m}n^{-1}),\\
V[\hat{\eta}_j(x_j)|\vec{X}_n]
&=&
\frac{1}{n}\vec{B}(x_j)^T\Gamma_j(\lambda_{jn})^{-1}\Gamma_j(0)\Gamma_j(\lambda_{jn})^{-1}\vec{B}(x_j)(1+o_P(1))\\
&=&O_P(K_n/n),\\
Cov(\hat{\eta}_i(x_i),\hat{\eta}_j(x_j))&=&O_P(n^{-1}),
\end{eqnarray*}
where $b_{j,a}(x_j)$ is given in Proposition \ref{app}, 
$$
b_{j,\lambda}(x)=-\frac{\lambda_{jn}}{n}\vec{B}(x)^T \Gamma_j(\lambda_{jn})^{-1}\Delta_m^\prime \Delta_m\vec{b}_{j0}=O\left(\frac{\lambda_n K_n^{1-m}}{n}\right).
$$
\end{theorem}

In Theorem \ref{mv}, the influence of $\gamma_n$ appears to be of only negligible order.
In actuality, we can use very small $\gamma_n$ as long as $H(\vec{b},\lambda_n,\gamma_n)$ is nonsingular.  
For example, Marx and Eilers (1998) used $\gamma_n=10^{-6}$. 
Thus, it is understood that the influence of $\gamma_n$ is small theoretically and numerically.
From Theorem \ref{mv}, the conditional Mean Squared Error(MSE) of $\hat{\eta}_j(x_j)$ can be obtained as follows. 

\begin{corollary}
Under the same assumption as Theorem \ref{mv}, it follows that
\begin{eqnarray*}
{\small
{\rm MSE}[\hat{\eta}_j(x_j)|\vec{X}_n]=E[\{\hat{\eta}_j(x_j)-\eta_j(x_j)\}^2|\vec{X}_n]=O_P\left(K_n^{-2(p+1)}+\lambda_{jn}^2K_n^{2(1-m)}n^{-2}\right)+O_P(K_nn^{-1}).
}
\end{eqnarray*}
Furthermore, the rate of convergence of the {\rm MSE} of $\hat{\eta}_j(x_j)$ becomes $O(n^{-(2p+2)/(2p+3)})$ by taking $K_n=O(n^{1/(2p+3)})$, $\lambda_{jn}=O(n^{\nu}), \nu\leq (p+m+1)/(2p+3)$.
\end{corollary}

Compared with the kernel estimator, the asymptotic order of MSE of the RCPS is the same as that of the local $p$th polynomial estimator when $p$ is odd and the number of knots in the spline methods and the bandwidth $h_n$ in the kernel methods are connected by $K_n/h_n^{-1}=O(1)$(see Opsomer (2000)). 
Lyapunov's condition of the central limit theorem yields the asymptotic normality of $[\hat{\eta}_1(x_1)\ \cdots\ \hat{\eta}_D(x_D)]^T$. 

\begin{theorem}\label{clt}
Suppose there exists $\delta\geq 2$ such that $E[|Y_i-c^\prime(\eta(\vec{x}_i))|^{2+\delta}|\vec{X}_i=\vec{x}_i]<\infty$. 
Furthermore, we assume $K_n=O(n^{1/(2p+3)})$ and $\lambda_{n}=O(n^\nu), \nu\leq (p+m+1)/(2p+3)$. 
Then, under the Assumptions, for any fixed point $\vec{x}=(x_1,\cdots,x_D)\in(0,1)^D$, as $n\rightarrow \infty$, 
\begin{eqnarray*}
\sqrt{\frac{n}{K_n}}
\left[
\begin{array}{c}
\hat{\eta}_{1,\gamma}(x_1)-\eta_1(x_1)-{\rm Bias}_1(x_1)\\
\vdots\\
\hat{\eta}_{D,\gamma}(x_D)-\eta_D(x_D)-{\rm Bias}_D(x_D)
\end{array}
\right]
\stackrel{d}{\longrightarrow}N_D\left( 
\vec{0}
,\Psi\right),
\end{eqnarray*}
where for $j=1,\cdots,D$, ${\rm Bias}_j(x_j)=b_{j,a}(x_j)+b_{j,\lambda}(x_j)$,
and $\Psi=\diag[\psi_1(x_1)\ \cdots\ \psi_D(x_D)]$ with
$$
\psi_j(x_j)=\lim_{n\rightarrow \infty}\frac{1}{K_n}\vec{B}(x_j)^T\Gamma_j(\lambda_{jn})^{-1}\Gamma_j(0)\Gamma_j(\lambda_{jn})^{-1}\vec{B}(x_j),\ \ j=1,\cdots,D.
$$
\end{theorem}

The proof of Theorem \ref{clt} is almost the same as that of Theorem 2 of Yoshida and Naito (2012). 
The asymptotic order of the bias and variance of the RCPS in Theorem \ref{mv} allows us to satisfy Lyapunov's condition for $[\hat{\eta}_1(x_1)\ \cdots\ \hat{\eta}_D(x_D)]^T$.

We note that an approximate pointwise confidence interval of $\eta_j(x_j)$ can be constructed by using the asymptotic distributional result of $\hat{\eta}_j(x_j)$. 
However, the asymptotic bias and variance of $\hat{\eta}_j(x_j)$ contain unknown variables and, hence, these should be estimated. 
For example, we replace $\vec{b}_0$ and $G_j$ with $\hat{\vec{b}}$ and $n^{-1}Z_j^T \hat{W}Z_j$, respectively, where $\hat{W}=\diag[c^{\prime\prime}(Z\hat{\vec{b}})]$. 
Furthermore, as it is the pilot estimator of the $(p+1)$th derivative of $\eta_j$, we can utilize the $(p+1)$th derivative of the RCPS $\hat{\eta}_j$ with $(p+2)$ or higher degree splines. 
Thus, we can construct the estimator $\widehat{{\rm Bias}}_j(x_j)$ and $\hat{\psi}_j(x_j)$ of ${\rm Bias}_j(x_j)$ and $\psi_j(x_j)$, respectively.
Consequently, we obtain an approximate confidence interval of $\eta_j(x_j)$ by the following Corollary.

\begin{corollary}\label{inter}
Under the same assumption as Theorem \ref{clt}, a $100(1-\alpha)\%$ asymptotic confidence interval of $\eta_j(x_j)$ at any fixed point $x_j\in(0,1)$ is
$$
\left[\hat{\eta}_j(x_j)-\widehat{{\rm Bias}}_j(x_j)- z_{\alpha/2}\sqrt{\hat{\psi}_j(x_j)},\  \hat{\eta}_j(x_j)-\widehat{{\rm Bias}}_j(x_j)+ z_{\alpha/2}\sqrt{\hat{\psi}_j(x_j)}\right],
$$
where $z_{\alpha/2}$ is the $(1-\alpha/2)$th normal percentile. 
\end{corollary}

\noindent{\bf Remark 1}
\ We see from the proof of Theorem \ref{mv} that the asymptotic form of $\hat{\eta}_j(x_j)$ can be written as 
\begin{eqnarray}
\hat{\eta}_j(x_j)-\eta_j(x_j)=\left\{\vec{B}(x_j)^T \Gamma_j(\lambda_{jn})^{-1}G_j(\vec{b}_0,\lambda_n,\gamma_n)+b_{j,a}(x_j)\right\}(1+o_P(1)) \label{asRCPS}
\end{eqnarray} 
under the same assumption as Theorem \ref{clt}, where $G_j(\vec{b}_0,\lambda_n,\gamma_n)$ is the $j$th $(K_n+p)$-subvector of $G(\vec{b}_0,\lambda_n,\gamma_n)$. 
From (2.8) of Kauermann et al. (2009), we see that $\hat{\eta}_j(x_j)$ and the penalized spline estimator based on the dataset $\{(y_i,x_{ij}):i=1,\cdots,n\}$ in GLM have the same asymptotic form. 
Thus, (\ref{asRCPS}) indicates that the asymptotic results of the RCPS in the GAM include those in the GLM. 
Note that in GLM($D=1$), we do not need to use the ridge penalty because the Hessian of the penalized log-likelihood of $\vec{b}$ is strictly convex.

\vspace{5mm}

\noindent{\bf Remark 2}
\ Claeskens et al. (2009) showed the asymptotic bias and variance of the penalized spline estimator in a regression model with $D=1$. 
They studied the asymptotics for penalized splines in the following two asymptotic scenarios: (a) the value $K_q$ appeared in their paper, less than 1, and (b) $K_q\geq 1$. 
In our setting, Assumption 5 guarantees case (a) and so Theorem \ref{mv} can be seen as the general version of Theorem 2 (a) of Claeskens et al. (2009) with respect to the model and dimension of covariates. 
If $\lambda_{jn}^{-1}$ is equal or smaller than the maximum eigenvalue of $(Z_j^TZ_j)^{-1}\Delta_m^T \Delta_m(Z_j^TZ_j)^{-1}$, the asymptotics for the penalized spline estimator in the GAM will be demonstrable, such as in Theorem 2 (b) of Claeskens et al. (2009).

\vspace{5mm}

\noindent{\bf Remark 3}
\ From Theorem \ref{clt}, it is understood that $[\hat{\eta}_1(x_1)\ \cdots\ \hat{\eta}_D(x_D)]^T$ are asymptotically mutually independent. 
Wand (1999) showed the asymptotic independence of the kernel estimator in additive models. 
Hence the penalized spline estimator and the kernel estimator have the same asymptotic property. 
The asymptotic independence of the joint distribution of $[\hat{\eta}_1(x_1)\ \cdots\ \hat{\eta}_D(x_D)]^T$ gives some justification for Corollary \ref{inter}, in which the approximate confidence interval is constructed based on the asymptotic result of the marginal distribution of $\hat{\eta}_j(x_j)$.

\vspace{5mm}

\noindent{\bf Remark 4}
\ Clearly, the penalized spline estimator can also be obtained via the backfitting algorithm. 
The asymptotic normality of the backfitting estimator can be shown, although it is not discussed in this paper. 
In additive models, Yoshida and Naito (2012) derived the asymptotic normality of the penalized spline estimator obtained by the backfitting algorithm.

\vspace{5mm}

\noindent{\bf Remark 5}
\ Theorems in this section have been shown for the RCPS with common $(p,K_n,m)$ in each covariate. 
When we construct $\hat{\eta}_j(x_j)$ using different $(p,K_n,m)$ in each $j$, the asymptotic normality of the RCPS can also be shown. 
In other words, for $\hat{\eta}_j(x_j)$ with $(p_j,K_{jn},m_j)$ which satisfy $(p_j,K_{jn},m_j)\not=(p_i,K_{in},m_i)$ ($j\not= i$), Theorems \ref{mv} and \ref{clt} hold.


\section{The mixed model representation}

In this section, we discuss the penalized spline estimator in relation to mixed models. 
We consider model (\ref{model1}) again with $\eta_j(x_j)$ approximated by a $p$th truncated spline model:
$$
\sum_{i=0}^{p} \beta_{i,j} x^{i} +\sum_{k=1}^{K_n-1} u_{k,j}(x-\kappa_k)^{p}_+,
$$ 
where $(x)_+=\max\{x,0\}$. 
We assume that the random vector $\vec{u}_j=(u_{1,j}\ \cdots\ u_{K_n-1,j})^T$ has density
$\vec{u}_j \sim N(\vec{0},\sigma_j^2 I)$ with $\sigma_j^2<\infty$,
and $\vec{u}_i$ and $\vec{u}_j$ are independent for $i\not= j$. 
Hence, $\vec{u}=[\vec{u}_1^T\ \cdots\ \vec{u}_D^T]^T$ distributes $N(\vec{0},\Sigma_u)$, where $\Sigma_u={\rm diag}[\sigma_1^2I\ \cdots\ \sigma_D^2I]$. 
Let 
\begin{eqnarray*}
X_j=
\left[ \begin{array}{cccc}
1 & x_{1j} & \cdots & x^p_{1j} \\
\vdots & \vdots & \ddots & \vdots \\
1 & x_{nj} & \cdots & x^p_{nj}
\end{array} \right], \quad
R_j=
\left[ \begin{array}{ccc}
(x_{1j}-\kappa_{1})^p_+ & \cdots & (x_{1j}-\kappa_{K_n-1})^p_+ \\ 
\vdots & \ddots & \vdots \\
(x_{nj}-\kappa_{1})^p_+ & \cdots & (x_{nj}-\kappa_{K_n-1})^p_+
\end{array} \right],
\end{eqnarray*}
$S_j=[X_j\ R_j]$, $S=[S_1\ \cdots\ S_D]$, $\vec{\beta}_j=(\beta_{0,j}\ \cdots\ \beta_{p,j})^T$, $\vec{\theta}_j=[\vec{\beta}_j^T\ \vec{u}_j^T]^T$ and $\vec{\theta}=[\vec{\theta}_1^T\ \cdots\ \vec{\theta}_D^T]^T$. 
The suggested joint density of $(\vec{y},\vec{u})$ can be written as 
\begin{eqnarray*}
f(\vec{y},\vec{u}:\vec{\beta})
&=&
f(\vec{y}|\vec{u}:\vec{\beta})f(\vec{u})\\
&=&
\exp\left[\vec{y}^T (S\vec{\theta})-\vec{1}^T c(S\vec{\theta})+\vec{1}^T h(\vec{y})\right]\frac{1}{\sqrt{(2\pi)^{D} |\Sigma_u|}}\exp\left[-\frac{1}{2}\vec{u}^T \Sigma_u^{-1}\vec{u}\right]\\
&=&
\exp\left[\vec{y}^T (S\vec{\theta})-\vec{1}^T c(S\vec{\theta})+\vec{1}^T h(\vec{y})\right]\frac{1}{\sqrt{(2\pi)^D |\Sigma_u|}}\exp\left[-\frac{1}{2}\vec{\theta}^T \Theta \vec{\theta}\right],
\end{eqnarray*}
where $\Theta={\rm diag}[\Theta_1\ \cdots\ \Theta_D]$ and $\Theta_j={\rm diag}[O_{p+1}\ \sigma_j^{-2}I]$.
As a convenient method of obtaining the estimator $\hat{\vec{\beta}}$ of $\vec{\beta}$ and the predictor $\hat{\vec{u}}$ of $\vec{u}$, the PQL is often used. 
In the PQL context, for given $\sigma_1^2,\cdots,\sigma_D^2$, $(\hat{\vec{\beta}},\hat{\vec{u}})$ is defined as the maximizer of 
\begin{eqnarray*}
\ell(\vec{\beta},\vec{u})=
 \frac{1}{n}\log f(\vec{y},\vec{u}|\vec{\beta})= 
\frac{1}{n}\{\vec{y}^T (S\vec{\theta})-\vec{1}^T c(S\vec{\theta})\}-\frac{1}{2n}\vec{\theta}^T \Theta\vec{\theta}+C(\Sigma_u),
\end{eqnarray*}
where $C(\Sigma_u)$ is not dependent on $\vec{\beta}$ and $\vec{u}$. 
Let $S(x)=(1\ x\ \cdots\ x^p\ (x-\kappa_1)_+^p\ \cdots\ (x-\kappa_{K_n-1})_+^p)^T$.  
Then, the PQL fit of $\eta_j(x_j)$ is defined as $\hat{\eta}_{j,P}(x_j)=S(x_j)^T \hat{\vec{\theta}}_j$, where $\hat{\vec{\theta}}_j=[\hat{\vec{\beta}}_j^T\ \hat{\vec{u}}_j^T]^T$.

We show the asymptotic distribution of $[\hat{\eta}_{1,P}(x_1)\ \cdots\ \hat{\eta}_{D,P}(x_D)]^T$. 
In order to achieve the asymptotic normality of the PQL fits, we consider the equivalence result between the $B$-spline model and the truncated spline model.
By the fundamental property of the $B$-spline function, there exists a $(K_n+p)$th invertible matrix $L_j$ such that $Z_j=S_jL_j$.  
Then we obtain $Z=SL$ and $S\vec{\theta}=Z\vec{b}$ where $L={\rm diag}[L_1\ \cdots\ L_D]$ and $\vec{b}=L^{-1}\vec{\theta}$. 
Furthermore, $\ell(\vec{\beta},\vec{u})$ can be rewritten as
\begin{eqnarray}
\ell(\vec{\beta},\vec{u})=\ell(\vec{b})
=
\frac{1}{n}\{\vec{y}^T (Z\vec{b})-\vec{1}^T c(Z\vec{b})\}-\frac{K_n^{2p}}{2n} \vec{b}^T Q_{p+1}(\Sigma_u) \vec{b}+C(\Sigma_u), \label{loglike}
\end{eqnarray}
where $Q_{p+1}(\Sigma_u)={\rm diag}[\sigma_1^{-2}\Delta_{p+1}^T\Delta_{p+1}\ \cdots\ \sigma_D^{-2}\Delta_{p+1}^T\Delta_{p+1}]$. 
Here we have used the fact that $\vec{\theta}^T\Theta\vec{\theta}=\vec{b}^T L^T\Theta L\vec{b}=K_n^{2p}\vec{b}^T Q_{p+1}(\Sigma_u)\vec{b}$. 
Claeskens et al. (2009) clarified the equality $\vec{\theta}_j^T \Theta_j\vec{\theta}_j=K_n^{2p}\vec{b}_j^T \Delta_{p+1}^T\Delta_{p+1}\vec{b}_j$.
By showing the asymptotic distribution of the maximizer $[\hat{\vec{b}}_{1,P}^T\ \cdots\ \hat{\vec{b}}_{D,P}^T]^T$ of $\ell(\vec{b})$, we obtain the asymptotic normality of $[\hat{\eta}_{1,P}(x_1)\ \cdots\ \hat{\eta}_{D,P}(x_D)]^T$, where 
\begin{eqnarray*}
\hat{\eta}_{j,P}(x_j)=S(\vec{x})^T \hat{\vec{\theta}}_j=S(\vec{x})^TL_j L_j^{-1} \hat{\vec{\theta}}_j=Z(\vec{x})^T \hat{\vec{b}}_{j,P}.
\end{eqnarray*}

\begin{theorem}\label{cltpql}
Suppose there exists $\delta\geq 2$ such that $E[|Y_i-c^{\prime}(\eta(\vec{x}_i))|^{2+\delta}|\vec{X}_i=\vec{x}_i]<\infty$ and $\eta_1,\cdots,\eta_D\in C^{p+1}$. 
Under $K_n=O(n^{1/(2p+3)})$ and $\sigma_j^{-2}=O(n^\nu),\nu<2/(2p+3)$,
for any fixed point $\vec{x}\in(0,1)^D$, as $n\rightarrow \infty$, 
\begin{eqnarray*}
\sqrt{\frac{n}{K_n}}
\left[
\begin{array}{c}
\hat{\eta}_{1,P}(x_1)-\eta_1(x_1)-{\rm Bias}_1(x_1)\\
\vdots\\
\hat{\eta}_{D,P}(x_D)-\eta_D(x_D)-{\rm Bias}_D(x_D)
\end{array}
\right]
\stackrel{d}{\longrightarrow}N_D\left( 
\vec{0}
,\Psi_P\right),
\end{eqnarray*}
where ${\rm Bias}_j(x_j)=b_{j,a}(x_j)+b_{j,\sigma}(x_j)$, $b_{j,a}(x_j)$ is that given in Proposition 1,
\begin{eqnarray*}
b_{j,\sigma}(x_j)=-\frac{K_n^{2p}}{n\sigma_j^{2}}\vec{B}(x_j)^T G_j^{-1}\Delta_{p+1}^{T}\Delta_{p+1}\vec{b}_{j0}=O_P(K_n^{p}/(n\sigma_j^2)),
\end{eqnarray*}
$\Psi_P=\diag[\psi_{1,P}(x_1)\ \cdots\ \psi_{D,P}(x_D)]$ and 
\begin{eqnarray*}
&&{\small \psi_{j,P}(x_j)}\\
&&{\small =\lim_{n\rightarrow \infty}\frac{1}{K_n}\vec{B}(x_j)^T (G_j+K_n^{2p}(n\sigma_j^2)^{-1}\Delta_{p+1}^T\Delta_{p+1})^{-1}G_j(G_j+K_n^{2p}(n\sigma_j^2)^{-1}\Delta_{p+1}^T\Delta_{p+1})^{-1}\vec{B}(x_j).
}
\end{eqnarray*}
\end{theorem}

\vspace{5mm}

\noindent{\bf Remark 6}
\ If we use $m=p+1$, the results of Theorem \ref{clt} are asymptotically equivalent to those of Theorem \ref{cltpql} by replacing $\lambda_{jn}$ with $K_n^{2p}\sigma_j^{-2}(j=1,\cdots,D)$.

\vspace{5mm}

\noindent{\bf Remark 7}
\ It should be noted that the maximum likelihood method or the restricted maximum likelihood method can be utilized for estimating $\sigma_j^2 (j=1,\cdots,D)$ by using pseudo data. 
These methods and estimation algorithm based on the Fisher-scoring algorithm are detailed by Breslow and Clayton (1993) and by Ruppert et al. (2003). 

\section{Applications}

We apply the approximate confidence interval of each covariate $\eta_j(x_j)$ to real datasets. 
In all examples, $(p,m)=(3,2)$ is adopted. 
The number of knots and the smoothing parameters are chosen via generalized cross-validation. 
As a pilot estimator of $\eta_j^{(4)}(x_j)$ in ${\rm Bias}_j(x_j)$, we utilize the 4th derivative of the RCPS with a 5th degree $B$-spline model.
To see the behavior of $\hat{\eta}_j(x_j)$, the partial residual plots 
$$
\hat{\eta}_j(x_{ij})+\hat{W}^{-1}(y_i-c^\prime(\hat{\eta}(\vec{x}_i))),
$$ 
for each $x_{ij}(j=1,\cdots,D)$ are displayed (see Cook and Dabrera (1998) and Landwehr et al. (1984)).

\subsection{Kyphosis data}

The kyphosis data set had 81 rows and 4 columns, representing data of children who have had corrective spinal surgery. 
This data is available in the software R (package ${\tt rpart}$).
For this data, the logistic model
\begin{eqnarray*}
Y_i&\sim&{\rm Bernoulli}\left(\frac{\exp[\eta_1(x_{i1})+\cdots+\eta_3(x_{i3})]}{1+\exp[\eta_1(x_{i1})+\cdots+\eta_3(x_{i3})]}\right),\ \ i=1,\cdots,81
\end{eqnarray*}
is assumed, where $Y_i$ is a factor with levels absent(0) or present(1) indicating whether a kyphosis was present (1) after the operation, $x_{i1}$ is the age in months, $x_{i2}$ is the number of vertebrae involved and $x_{i3}$ is the number of the first (topmost) vertebra operated on. 
We construct the RCPS with $\gamma_n=10^{-6}$ and the approximate confidence intervals for each $\eta_j(x_j)$.  

\begin{figure}
\begin{center}
\includegraphics[width=45mm,height=45mm]{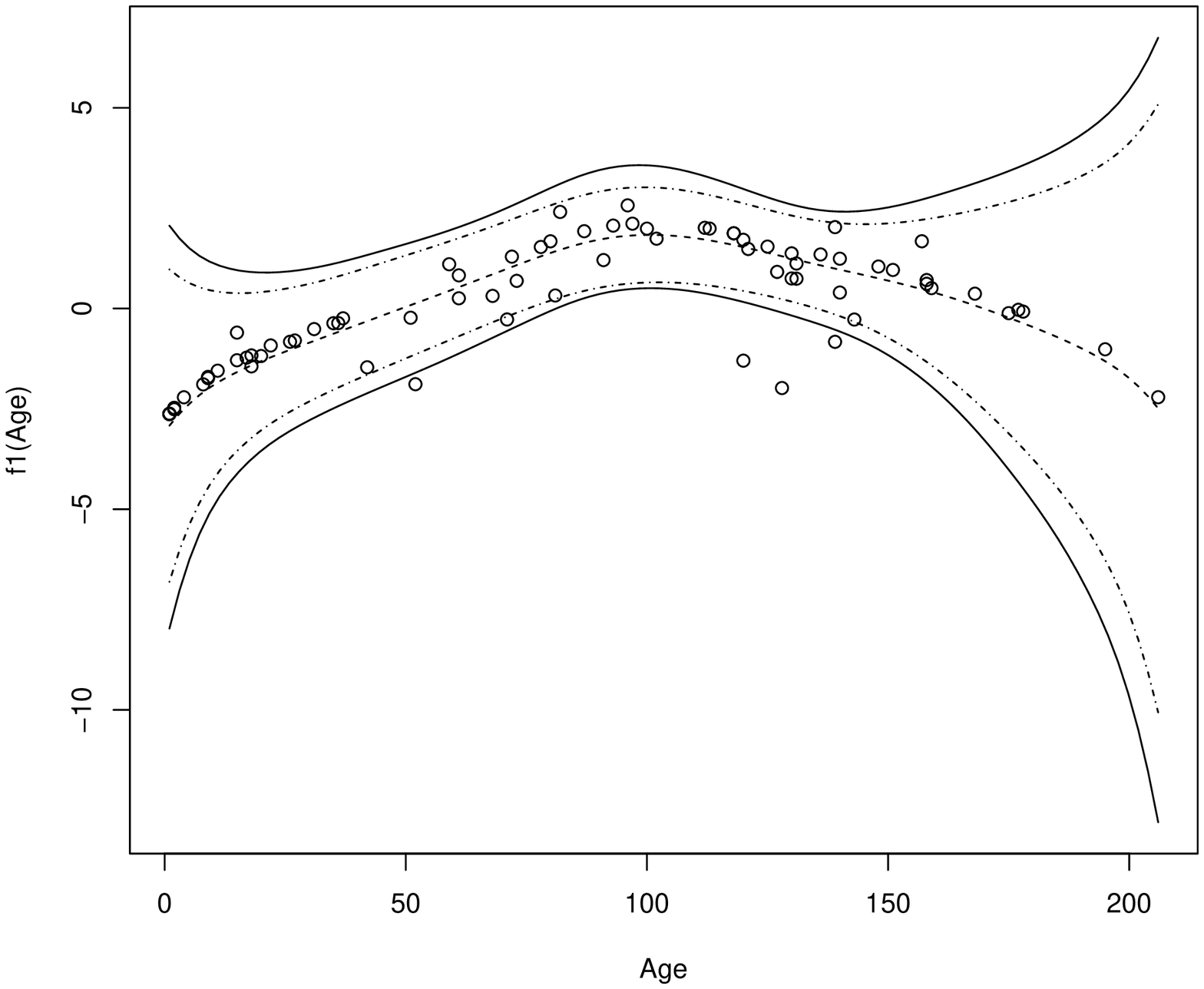}
\includegraphics[width=45mm,height=45mm]{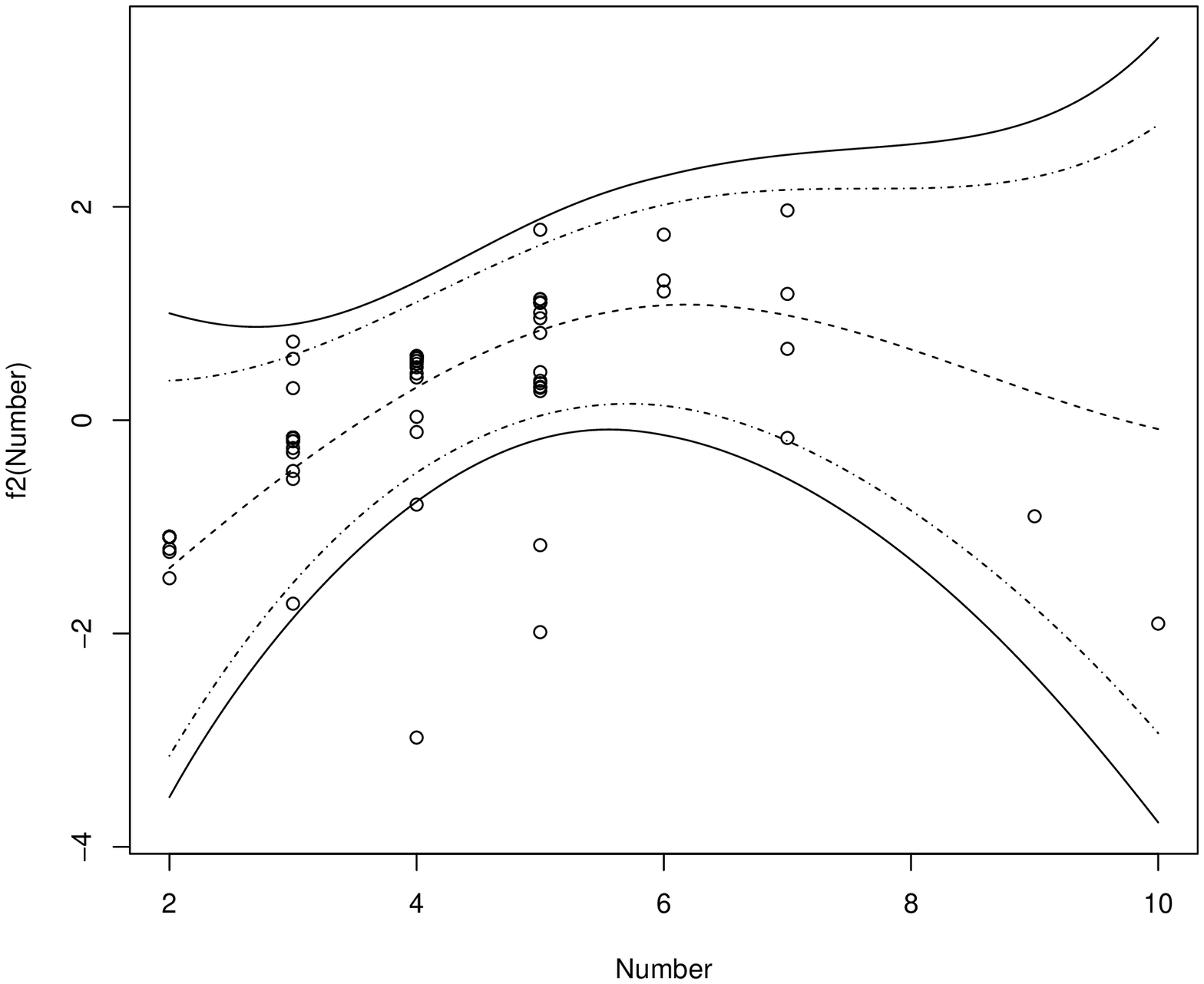}
\includegraphics[width=45mm,height=45mm]{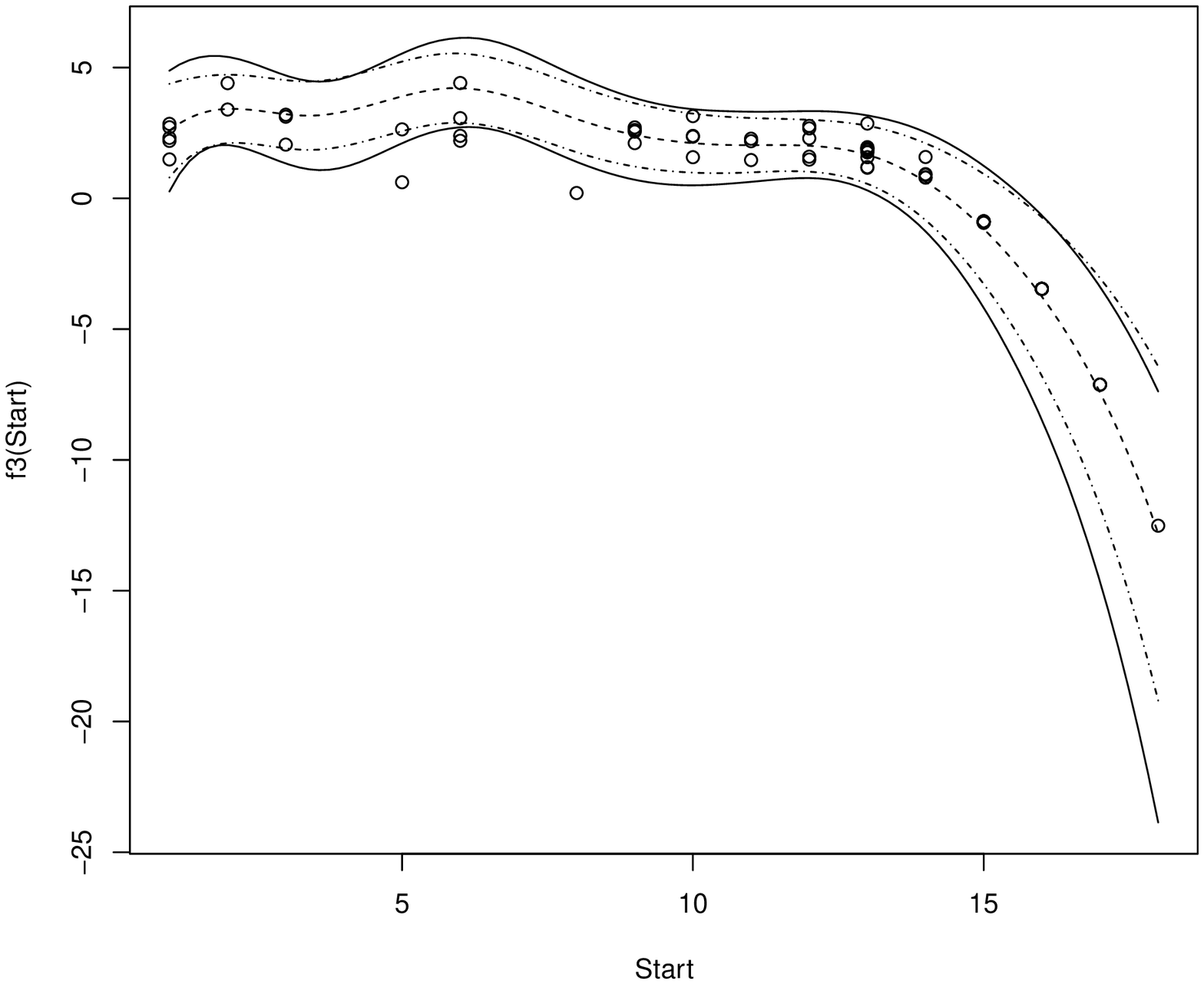}
\end{center}
\caption{Plots of the kyphosis data with the RCPS (dashed), the 99$\%$ approximate confidence interval(solid), $\hat{\eta}_j\pm 2\times$(standard error) (dot-dashed line) and the partial residuals. 
The left, middle and right panels are for $\eta_1(x_1)$, $\eta_2(x_2)$ and $\eta_3(x_3)$, respectively. \label{kyphosis}}
\end{figure}

In Fig. \ref{kyphosis}, for $j=1,2,3$, the RCPS $\hat{\eta}_j(x_j)$, the 99$\%$ approximate pointwise confidence interval 
\begin{eqnarray*}
\left[\hat{\eta}_j(x_j)-\widehat{{\rm Bias}}_j(x_j)- 2.58\sqrt{\hat{\psi}_j(x_j)},\  \hat{\eta}_j(x_j)-\widehat{{\rm Bias}}_j(x_j)+ 2.58\sqrt{\hat{\psi}_j(x_j)}\right],
\end{eqnarray*}
and the partial residual are all illustrated. 
For comparison, $\eta_j\pm 2\times$(standard error): 
\begin{eqnarray*}
\left[\hat{\eta}_j(x_j)- 2\sqrt{\hat{\psi}_j(x_j)},\  \hat{\eta}_j(x_j)+ 2\sqrt{\hat{\psi}_j(x_j)}\right],\ j=1,2,3
\end{eqnarray*}
are also superimposed. 
In all covariates, smooth intervals are obtained. 
Marx and Eilers (1998) illustrated the RCPS and $\eta_j\pm 2\times$(standard error) for the same dataset in Fig.4 of their paper. 
Our results and theirs are similar. 
However, our interval is wiggles a bit because the asymptotic bias is corrected in each covariate.

\subsection{Air Pollution and Mortality data}

\begin{figure}
\begin{flushleft}
\includegraphics[width=45mm,height=40mm]{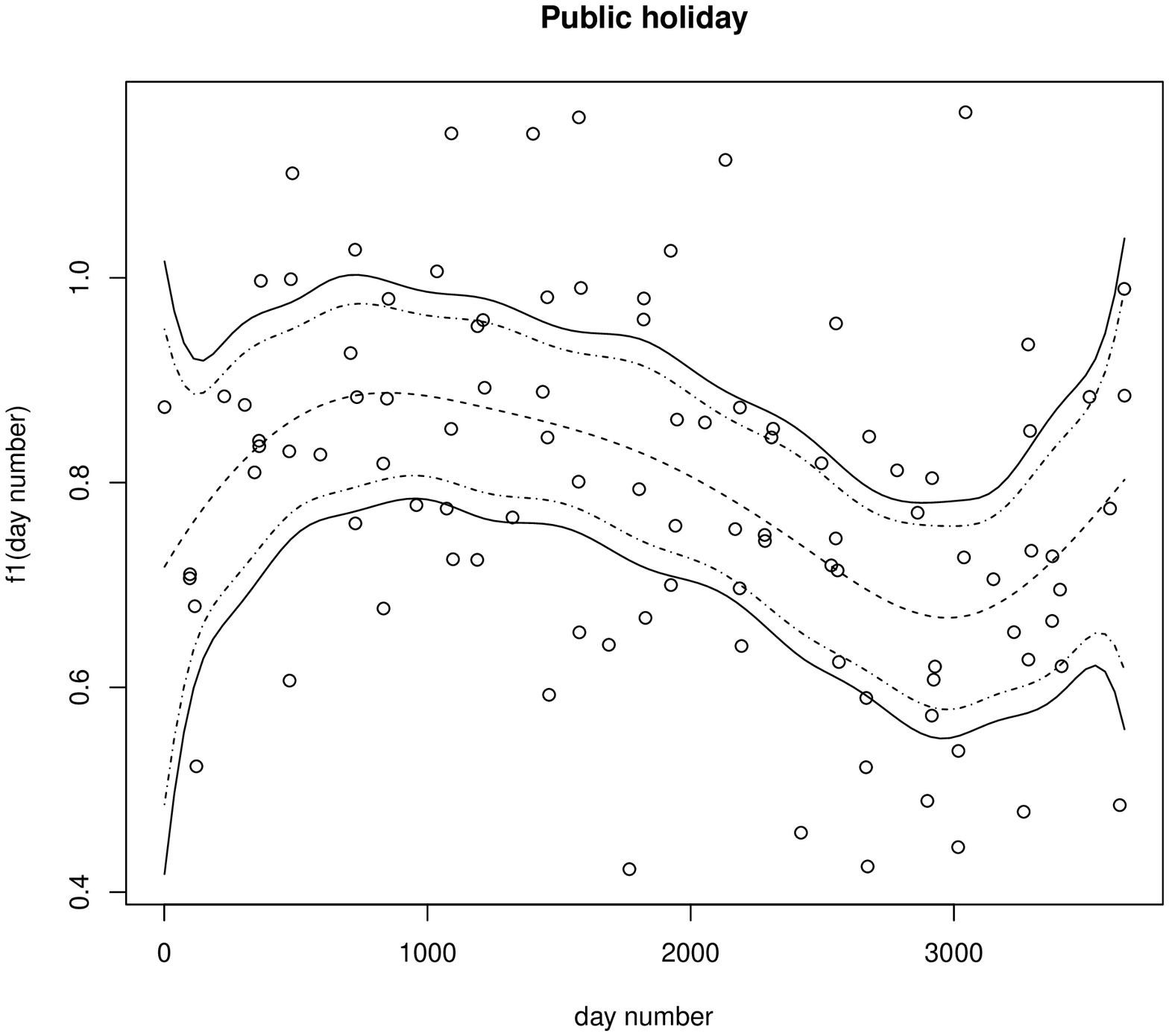}
\includegraphics[width=45mm,height=40mm]{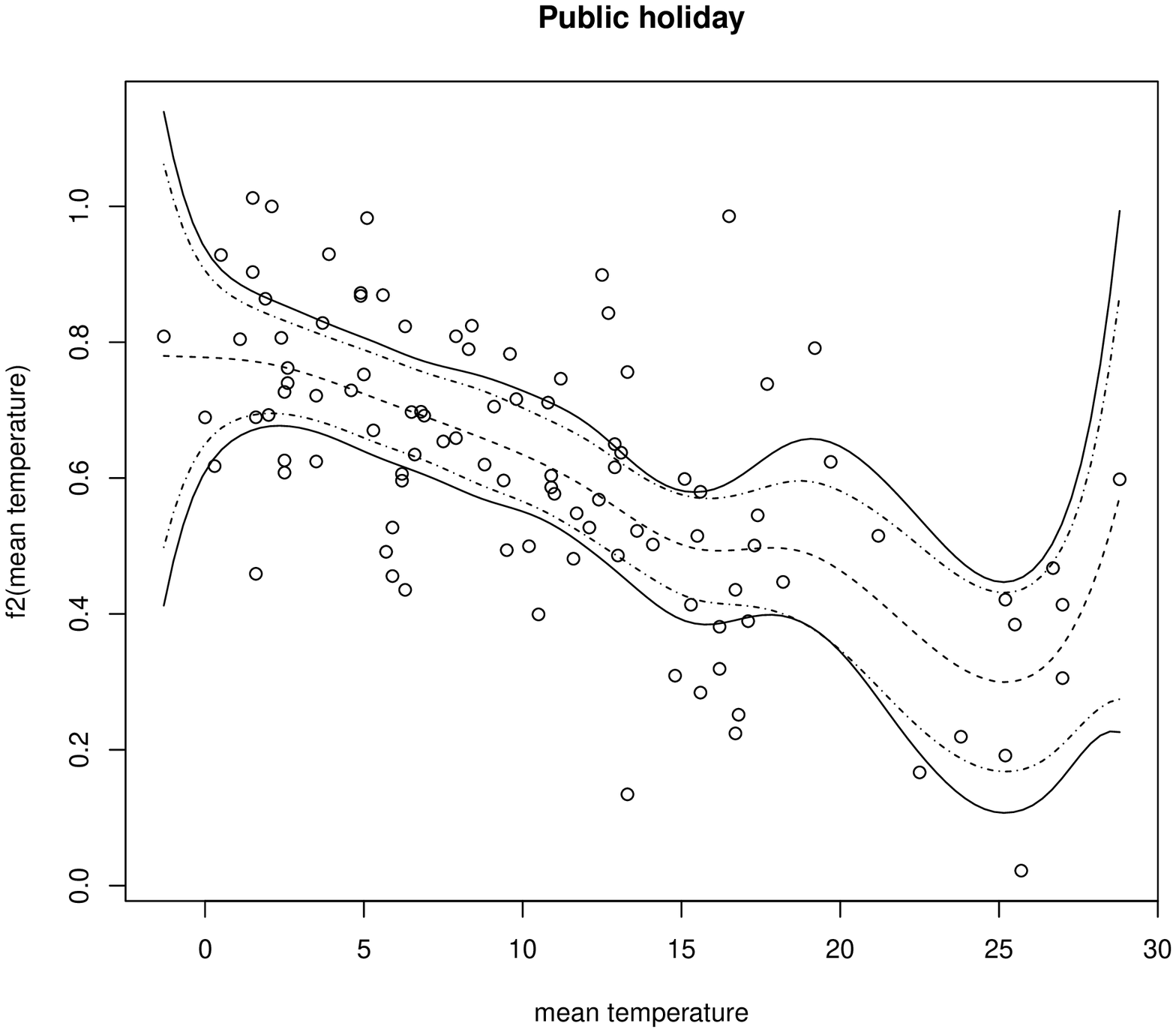}
\includegraphics[width=45mm,height=40mm]{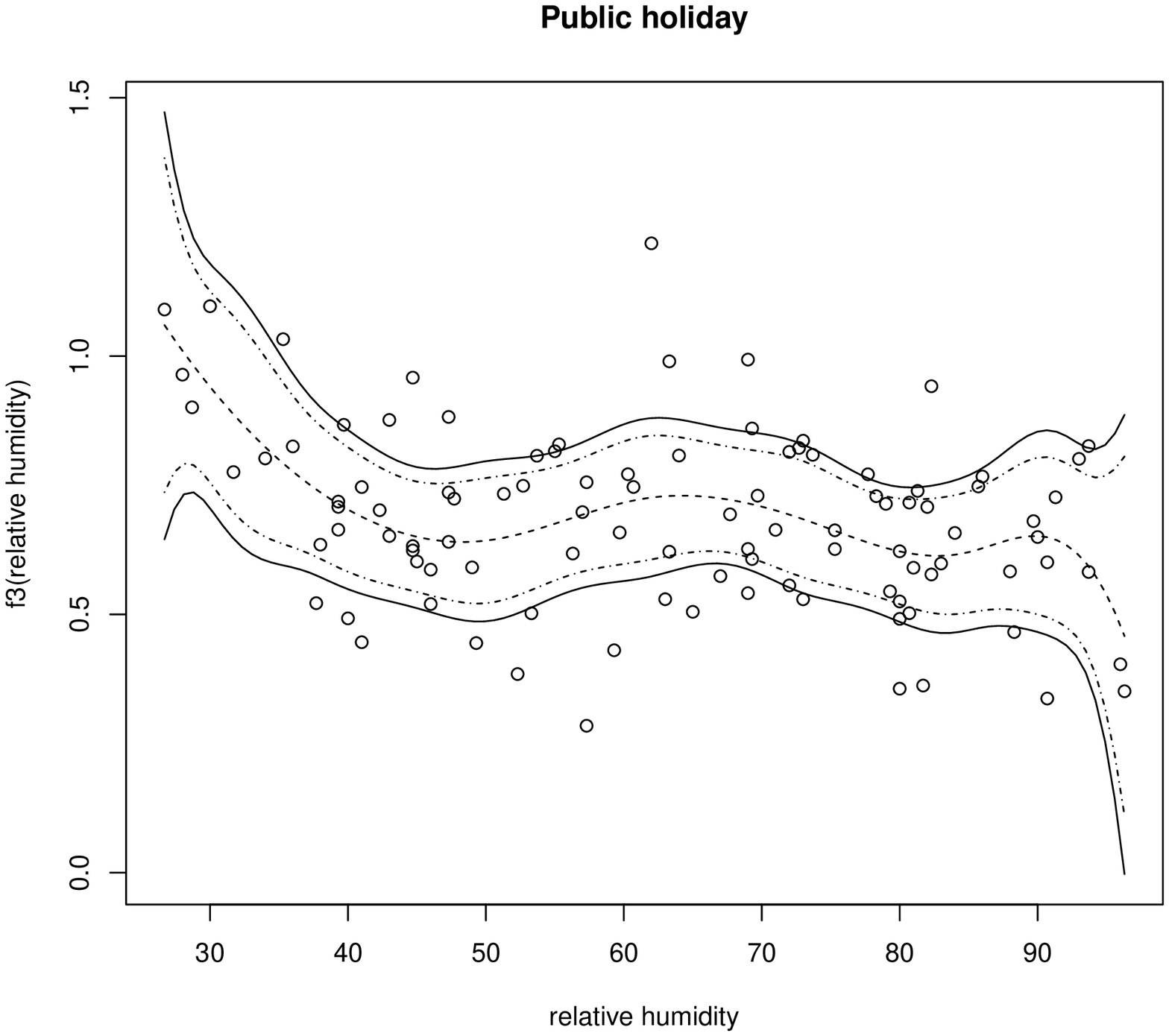}\\
\includegraphics[width=45mm,height=40mm]{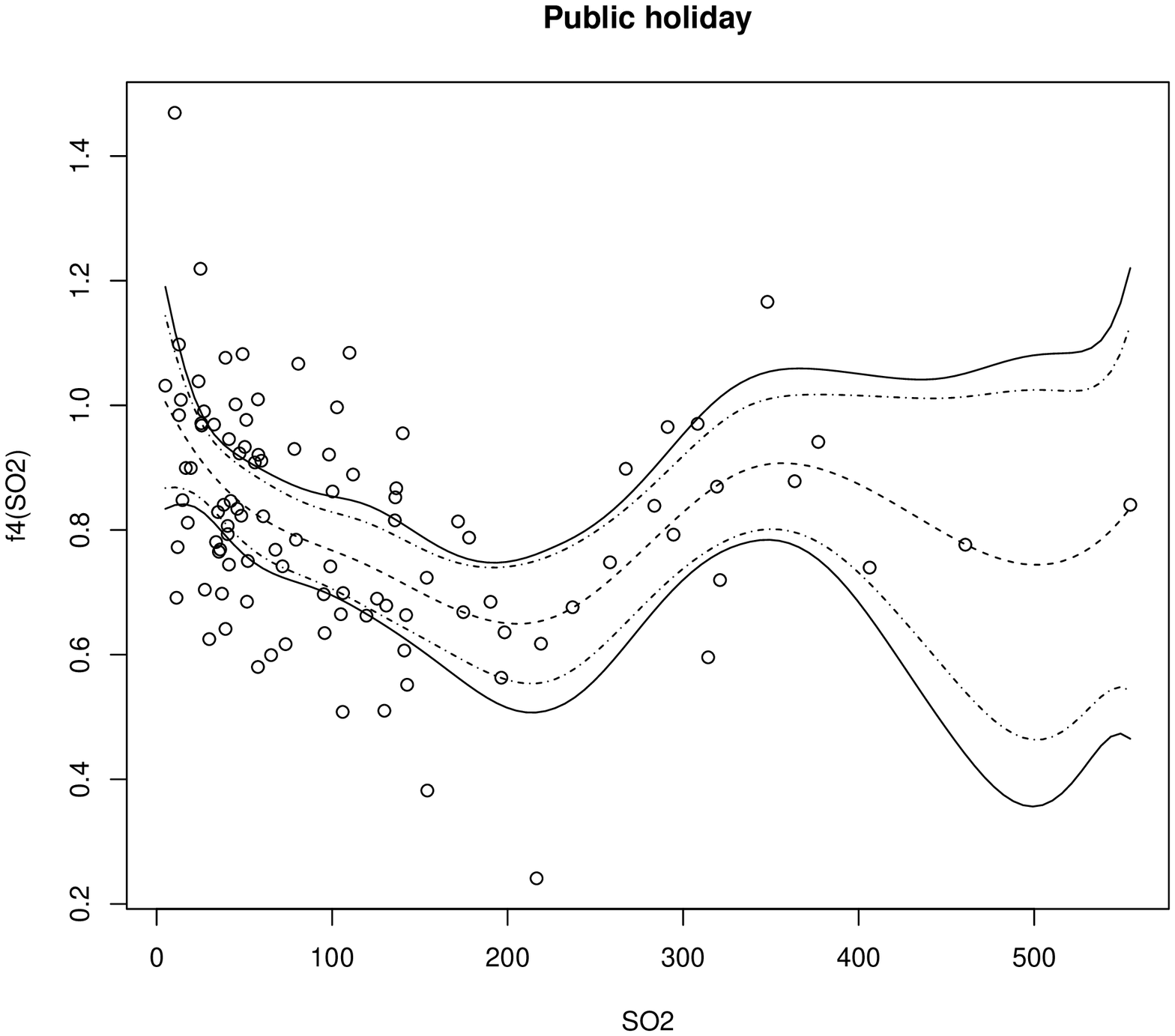}
\includegraphics[width=45mm,height=40mm]{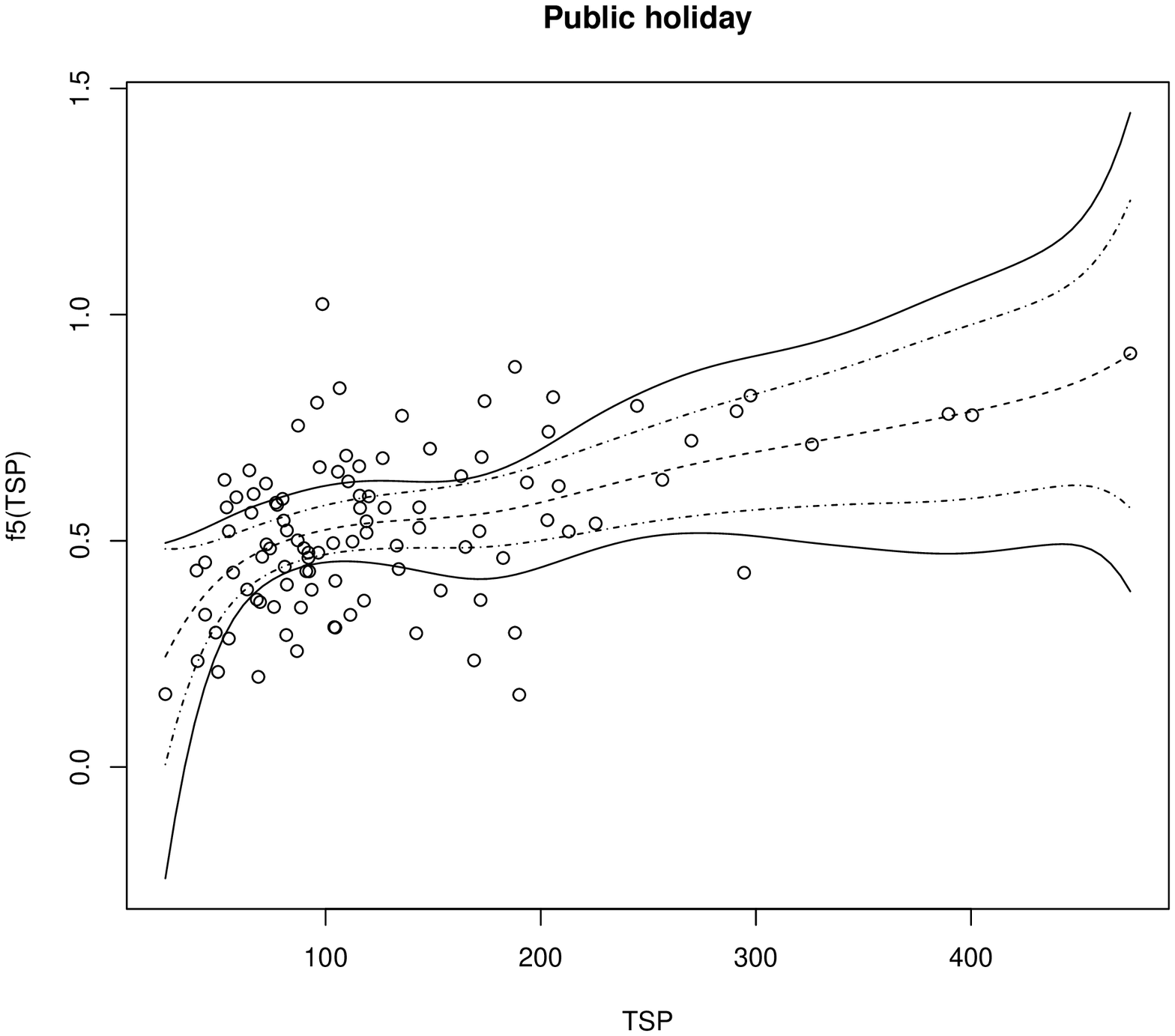}
\end{flushleft}
\caption{Plots of air pollution and mortality data with the RCPS, the 99$\%$ approximate confidence interval, $\hat{\eta}_j\pm 2\times$(standard error) and the partial residuals.  
 \label{air}}
\end{figure}

This data set contained air pollution and mortality data for the city of Milan, Italy, over 3652 consecutive days (i.e., 10 consecutive years: 1st January, 1980 to 30th December, 1989). 
The original data is available on the web site of Ruppert et al. (2003). 
The relationship between the number of deaths in a day and some explanatory variables is modeled as follows 
\begin{eqnarray*}
Y_i \sim {\rm Poisson}[\exp(\eta_1(x_{i1})+\cdots+\eta_5(x_{i5}))],\ \ i=1,\cdots,102,
\end{eqnarray*}
where $Y_i$ is the total number of deaths in a day, $x_{i1}$ is the number of days since 31st December, 1979, $x_{i2}$ the mean daily temperature in degrees celcius, $x_{i3}$ is the relative humidity, $x_{i4}$ is a measure of sulfur dioxide levels (SO2) in ambient air and $x_{i5}$ is the total amount of suspended particles (TSP) in ambient air. 
All of these have been measured on public holidays within the 3652 days, giving a sample size of $n=102$. 
We constructed the RCPS of $\eta_j(x_j)$ and the 99$\%$ approximate confidence intervals. 
In Fig. \ref{air}, the RCPS, the $99\%$ approximate confidence intervals, $\hat{\eta}_j\pm2\times$(standard error) and the partial residual for each $x_j$ are illustrated. 
We see that the effect of $\widehat{{\rm Bias}}_j(x_j)$ is somewhat large for all covariates.

\section{Simulation}

In this section, we validate Theorem \ref{clt} numerically by simulation. 
The true natural parameter utilized in the simulation is defined as $\eta(\vec{x})=\eta_1(x_1)+\eta_2(x_2)+\eta_3(x_3)$, where $\eta_1(x_1)=\sin(2\pi x_1)$, $\eta_2(x_2)=2\cos(2\pi x_2)$ and $\eta_3(x_3)=\sin^2((\pi/2) x_3)$. 
The design points $(x_{i1},x_{i2},x_{i3})$ are created by 
\begin{eqnarray*}
\left[
\begin{array}{c}
x_{i1}\\
x_{i2}\\
x_{i3}
\end{array}
\right]
=
\left[
\begin{array}{ccc}
(1+\rho+\rho^2)^{-1}&0&0\\
0&(1+2\rho)^{-1}&0\\
0&0&(1+\rho+\rho^2)^{-1}
\end{array}
\right]
\left[
\begin{array}{ccc}
1&\rho&\rho^2\\
\rho&1&\rho\\
\rho^2&\rho&1
\end{array}
\right]
\left[
\begin{array}{c}
z_{i1}\\
z_{i2}\\
z_{i3}
\end{array}
\right],
\end{eqnarray*} 
where $z_{ij}(i=1,\cdots,n,j=1,2,3)$ are generated independently from $U(0,1)$, the uniform distribution on $[0,1]$. 
We prepared two types of the design, with (i) $\rho=0$ and (ii) $\rho=0.2$. 
Then, the true functions are corrected to satisfy $E[\eta_j(X_j)]=0$ in each (i) and (ii).  
The response $Y_i$ is generated from 
\begin{eqnarray}
Y_i\sim{\rm Bernoulli}\left(\frac{\exp[\eta_1(x_{i1})+\eta_2(x_{i2})+\eta_3(x_{i3})]}{1+\exp[\eta_1(x_{i1})+\eta_2(x_{i2})+\eta_3(x_{i3})]}\right),\ \ i=1,\cdots,n. \label{ber1}
\end{eqnarray}

Our purpose is to compare the density of $N(0,1)$ and the kernel density estimate of the simulated $U_j (j=1,2,3)$, as well as the density of $N_2(\vec{0},I_2)$ and the kernel density estimate of the simulated $[U_1,U_2]^T$, $[U_1,U_3]^T$ and $[U_2,U_3]^T$ to validate Theorem \ref{clt}, where
\begin{eqnarray}
\left[
\begin{array}{c}
U_1\\
U_2\\
U_3
\end{array}
\right]
=
\sqrt{\frac{n}{K_n}}
\left[
\begin{array}{c}
\displaystyle\frac{\hat{\eta}_1(x_1)-\eta_1(x_1)-\widehat{{\rm Bias}}_1(x_1)}{\hat{\psi}_1(x_1)}\\
\displaystyle\frac{\hat{\eta}_2(x_2)-\eta_2(x_2)-\widehat{{\rm Bias}}_2(x_2)}{\hat{\psi}_2(x_2)}\\
\displaystyle\frac{\hat{\eta}_3(x_3)-\eta_3(x_3)-\widehat{{\rm Bias}}_3(x_3)}{\hat{\psi}_3(x_3)}
\end{array}
\right]. \label{2norm}
\end{eqnarray}
Here, 
\begin{eqnarray*}
\hat{\psi}_j(x_j)=\frac{1}{K_n}\vec{B}(x_j)^T\hat{\Gamma}_j(\lambda_{jn})^{-1}\hat{\Gamma}_j(0)\hat{\Gamma}_j(\lambda_{jn})^{-1}\vec{B}(x_j),
\end{eqnarray*}
$\hat{\Gamma}_j(\lambda_{jn})=n^{-1}(Z_j^T\hat{W}Z_j+\lambda_{jn}\Delta_m^\prime \Delta_m)$. 
For $j=1,2,3$, $\widehat{{\rm Bias}}_j(x_j)$ is constructed using the same method as that in the previous section. 
The bandwidth discussed by Sheather and Jones (1991) is utilized for kernel density estimates. 
The simulation algorithm described as follows: 

\vspace{5mm}

{\small

\begin{enumerate}
\item[Step 1]\ For $j=1,2,3$ and $i=1,\cdots,n$, generate $x_{ij}$ from (i) or (ii).
\item[Step 2]\ Generate the data $\{(y_{i},\vec{x}_{i})|i=1,\cdots,n\}$ from (\ref{ber1}).
\item[Step 3]\ Calculate $\hat{\eta}_j(x_j) (j=1,2,3)$ at a fixed point $(x_1,x_2,x_3)=(0.5,0.5,0.5)$.
\item[Step 4]\ Calculate the values of (\ref{2norm}).
\item[Step 5]\ Iterate from Step 2 to Step 4, 10000 times.
\item[Step 6]\ Draw the kernel density estimate of $U_1, U_2$ and $U_3$ and compare with the density of $N(0,1)$. 
\item[Step 7]\ Draw the kernel density estimate of $[U_1\ U_2]^T$, $[U_1\ U_3]^T$ and $[U_2\ U_3]^T$, and compare with the density of $N_2(\vec{0},I_2)$.
\end{enumerate}
}

\vspace{5mm}

To construct $\hat{\eta}_j(x_j) (j=1,2,3)$, we utilize the cubic $B$-spline ($p=3$) and the second difference matrix ($m=2$). 
Furthermore $K_n=2\lceil n^{2/5}\rceil$, $\lambda_{1n}=0.1\sqrt{n/K_n}$, $\lambda_{2n}=0.01\sqrt{n/K_n}$ and $\lambda_{3n}=\sqrt{n/K_n}$ are used. 
The ridge parameter is chosen as $\gamma_n=10^{-4}$. 
The sample sizes are set at $n=100$ and $n=1000$. 

\begin{figure}
\begin{center}
\includegraphics[width=45mm,height=45mm]{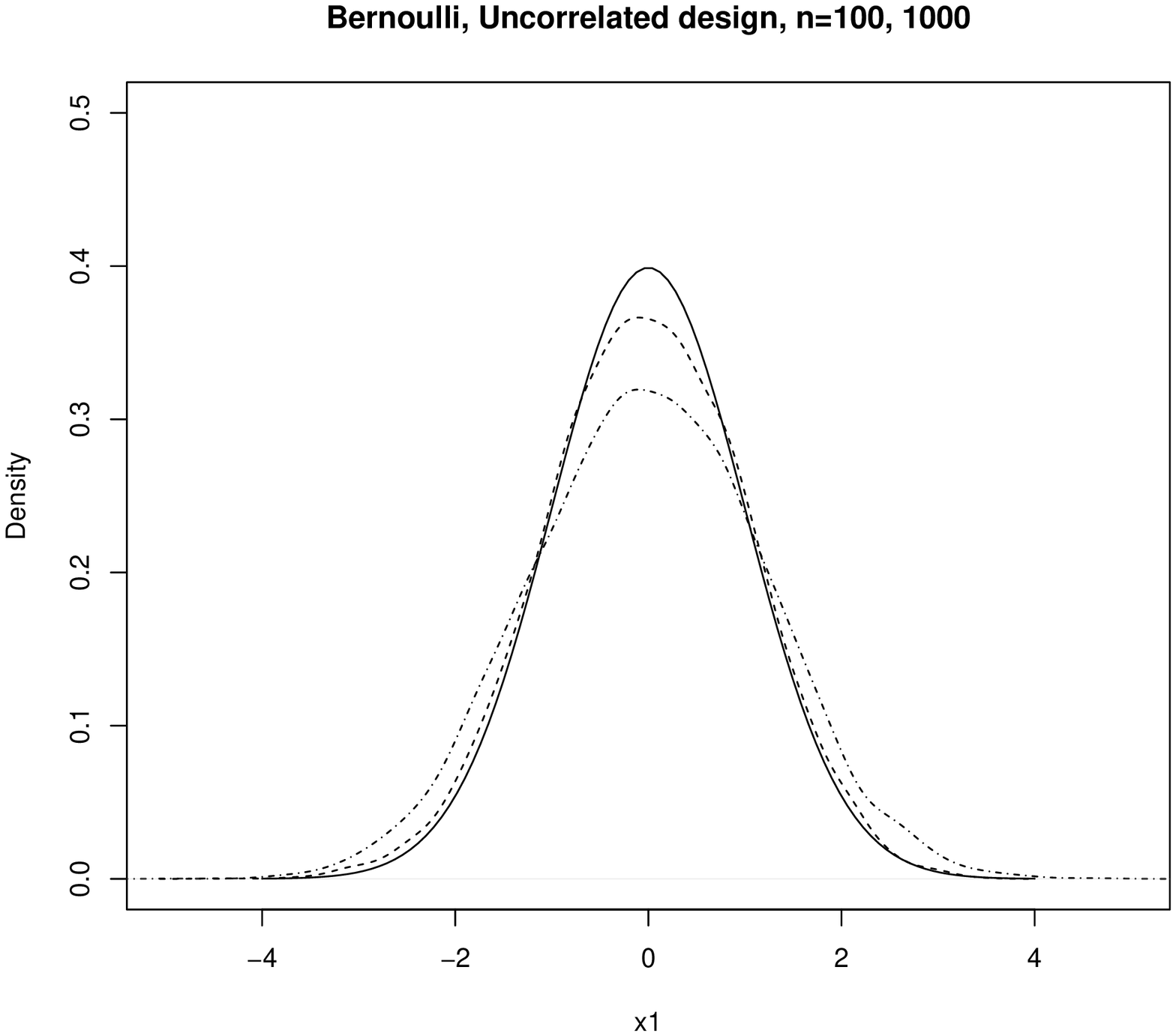}
\includegraphics[width=45mm,height=45mm]{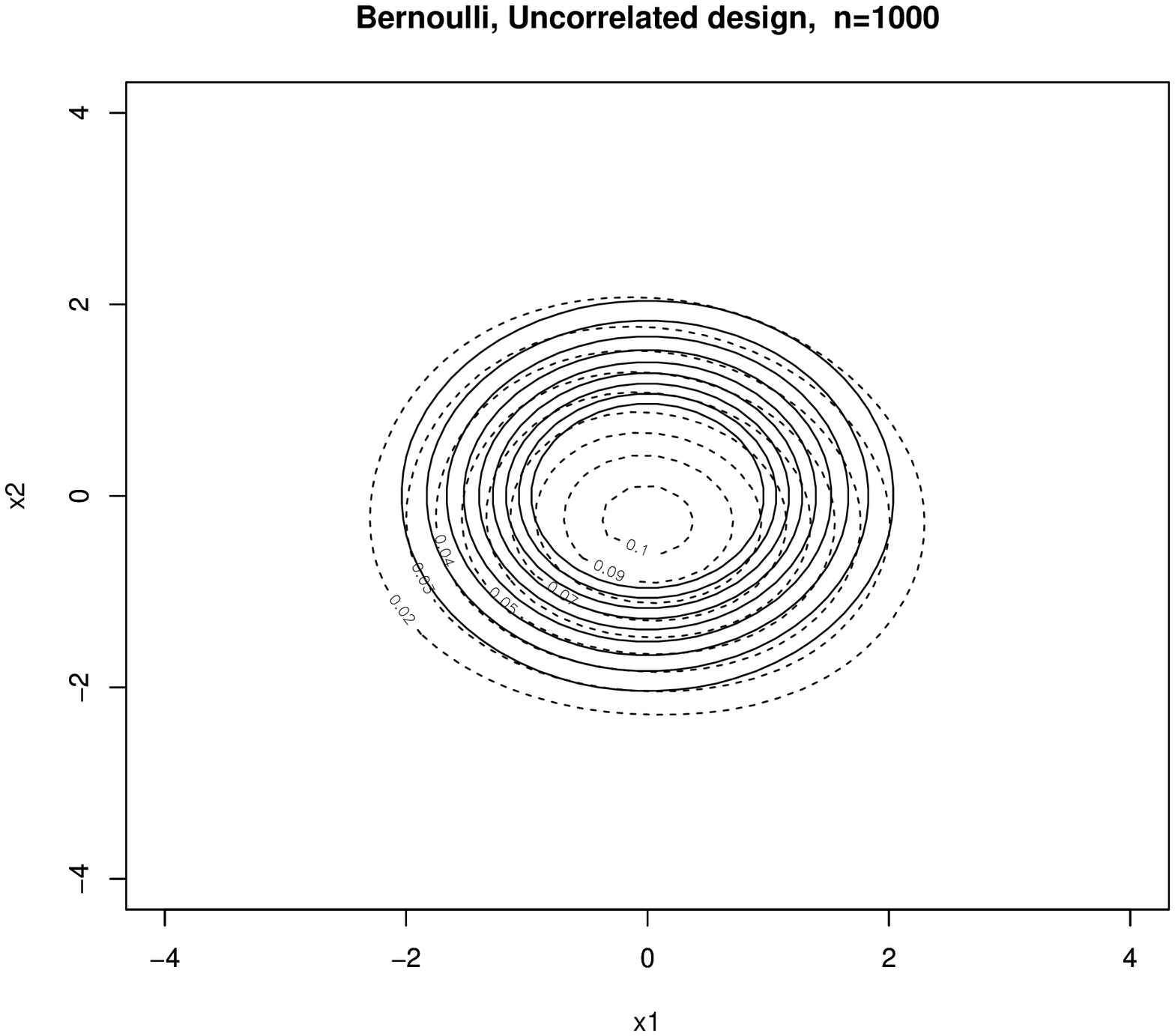}
\includegraphics[width=45mm,height=45mm]{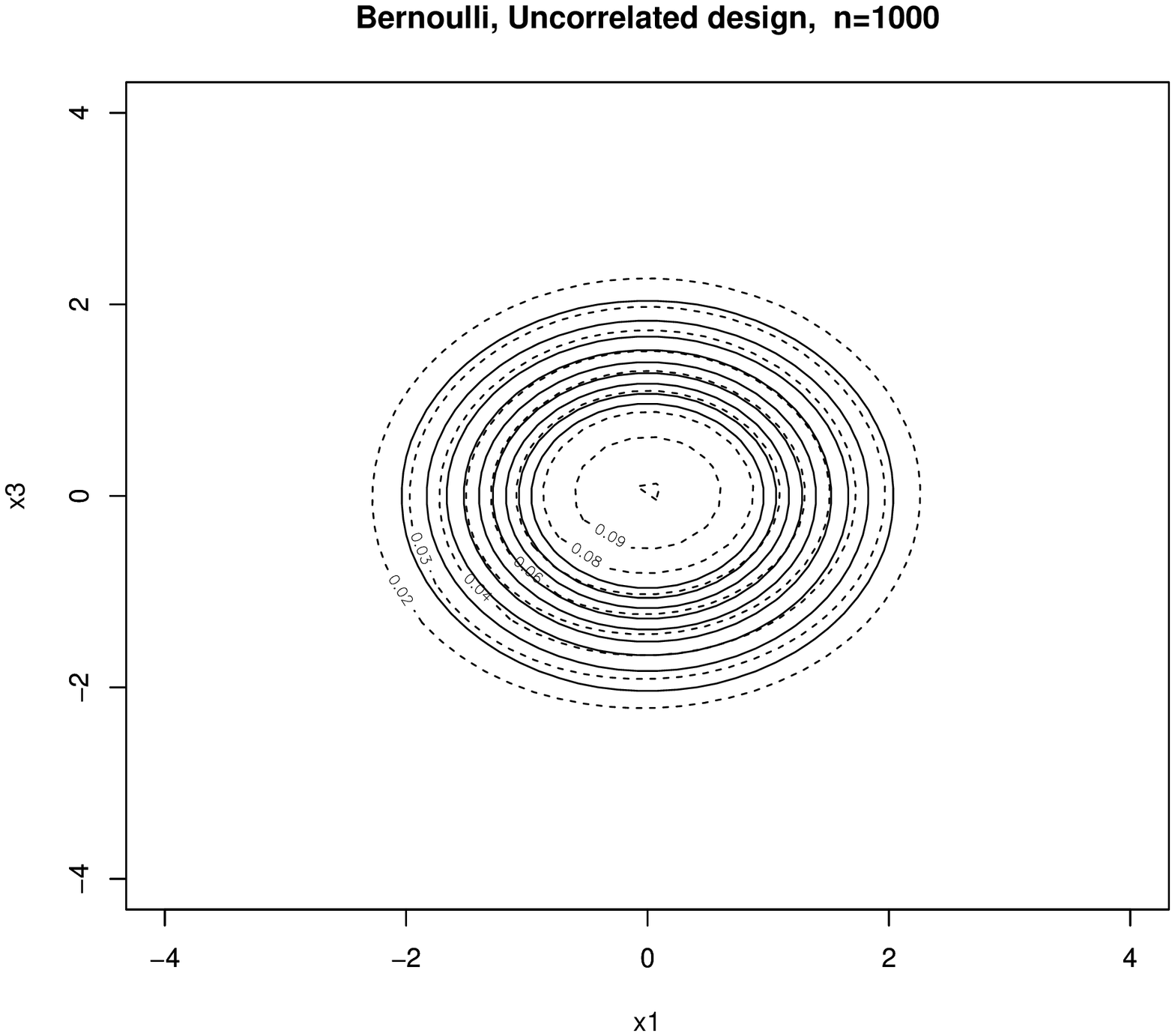}\\
\includegraphics[width=45mm,height=45mm]{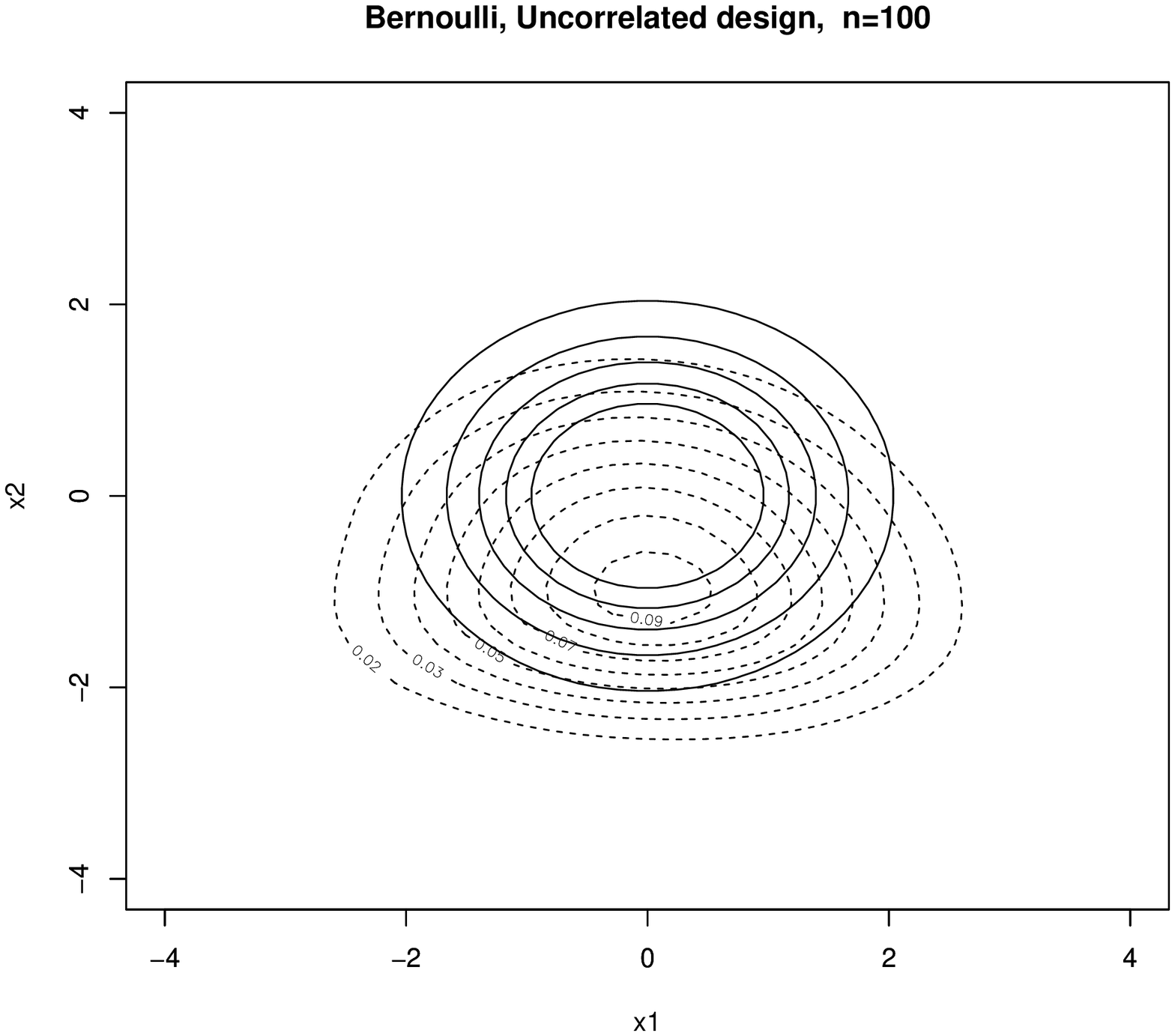}
\includegraphics[width=45mm,height=45mm]{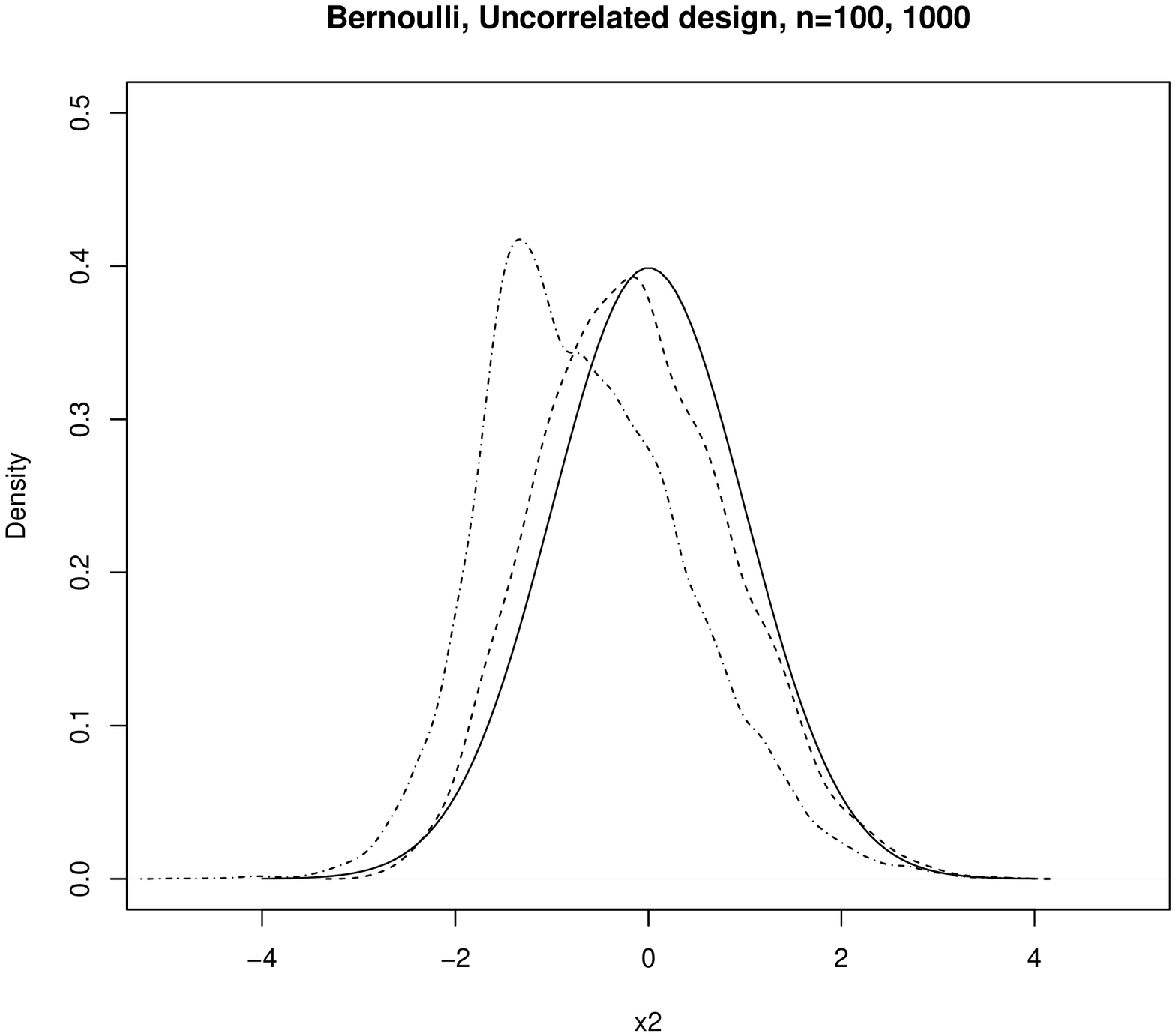}
\includegraphics[width=45mm,height=45mm]{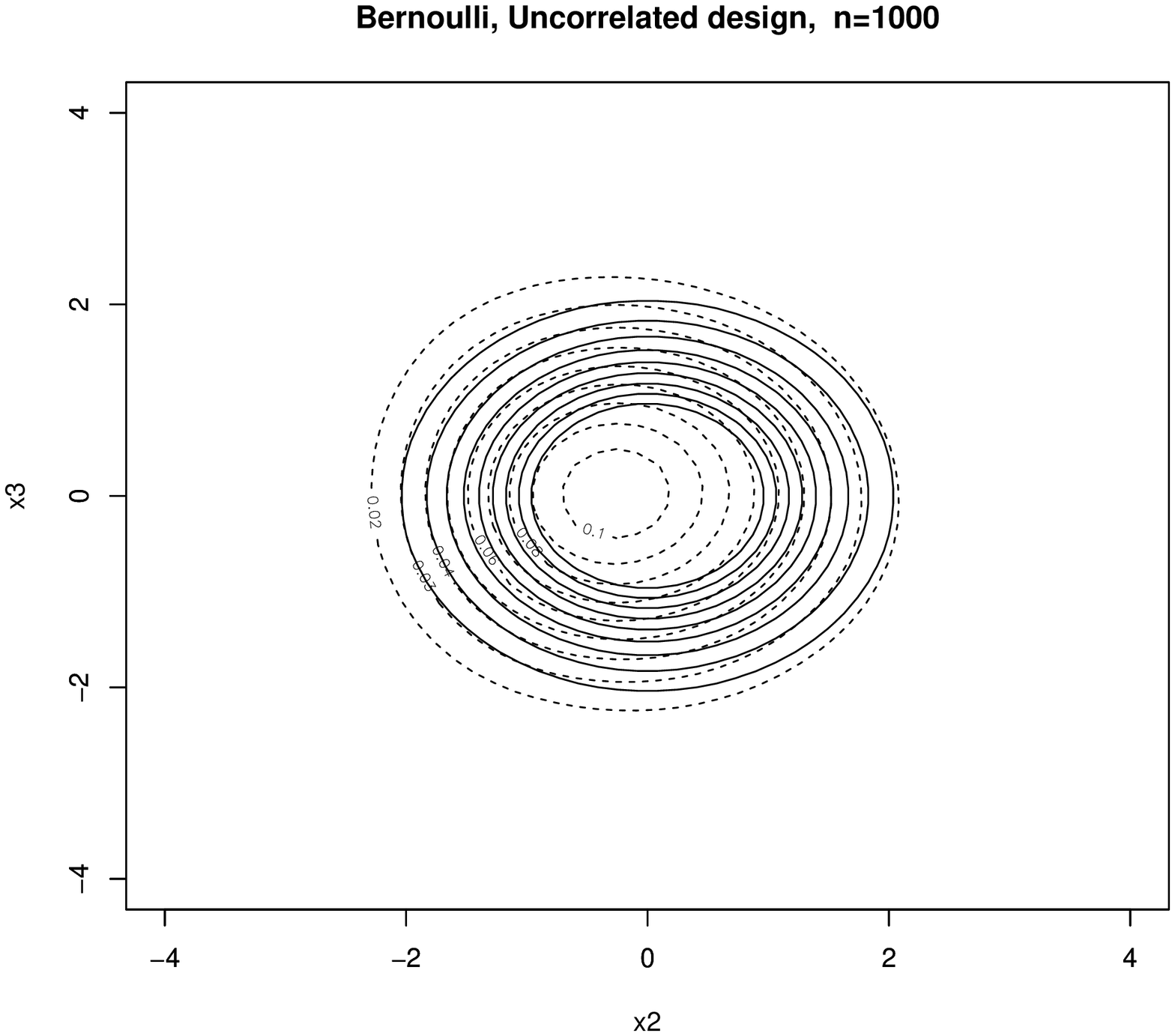}\\
\includegraphics[width=45mm,height=45mm]{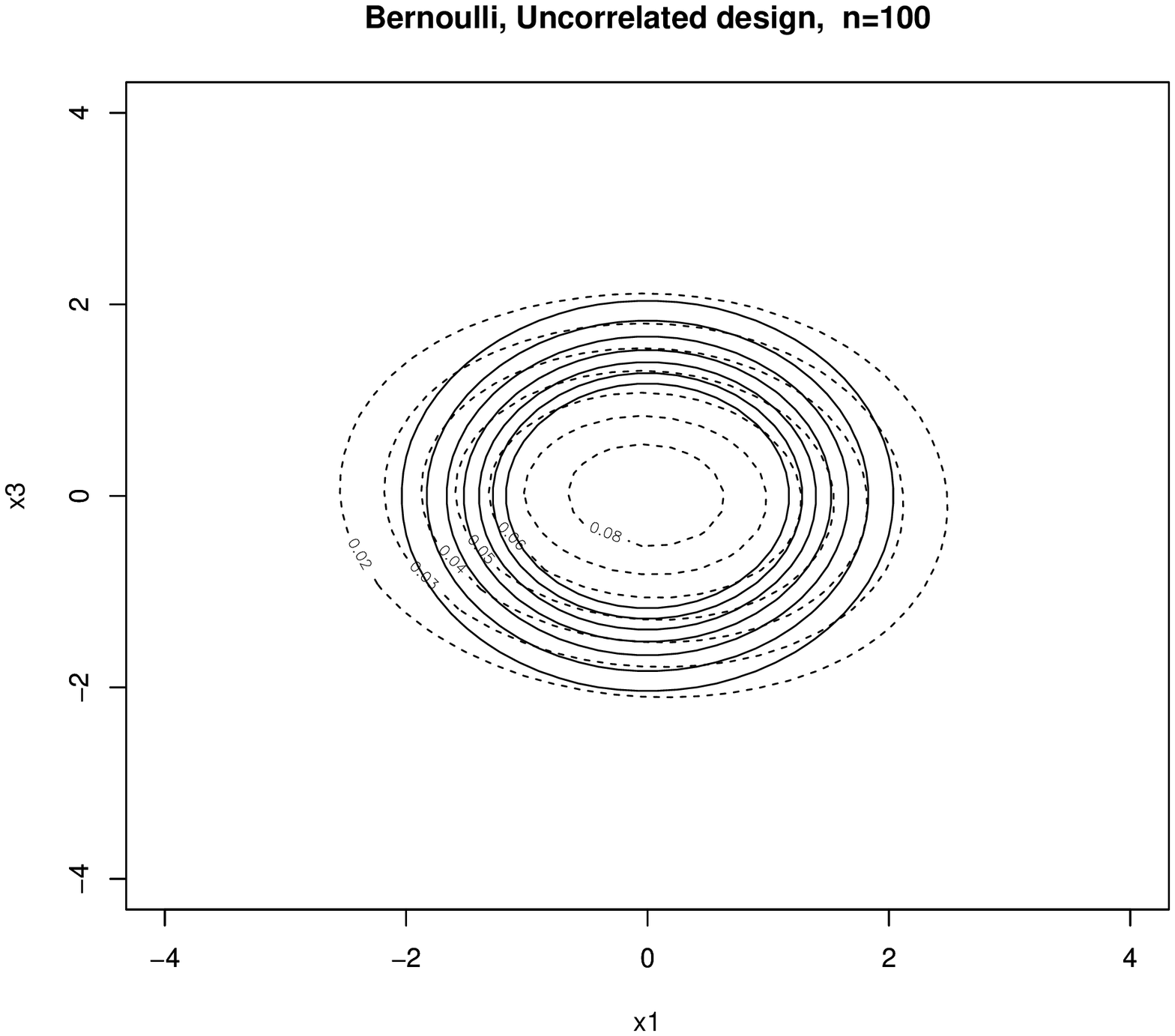}
\includegraphics[width=45mm,height=45mm]{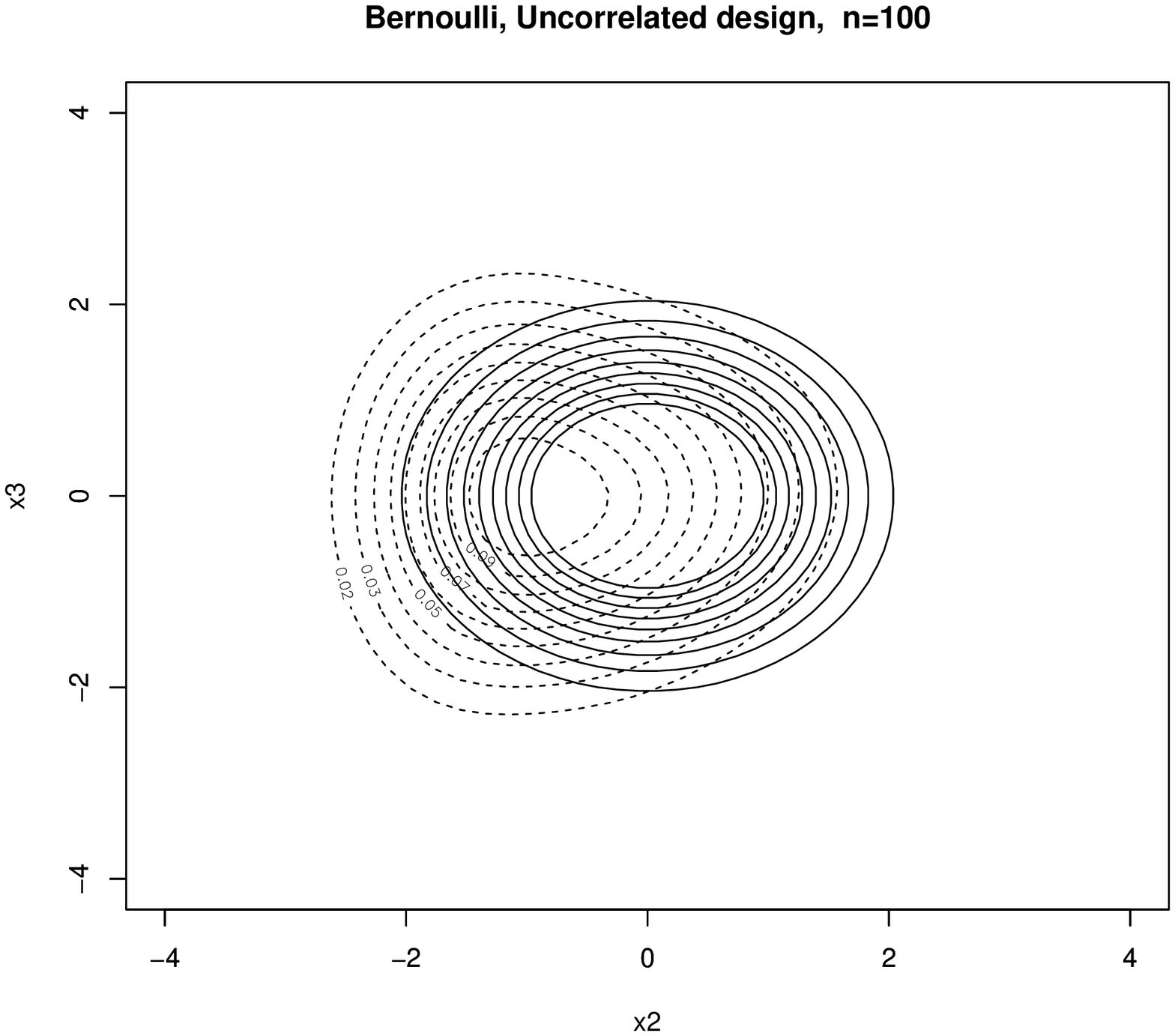}
\includegraphics[width=45mm,height=45mm]{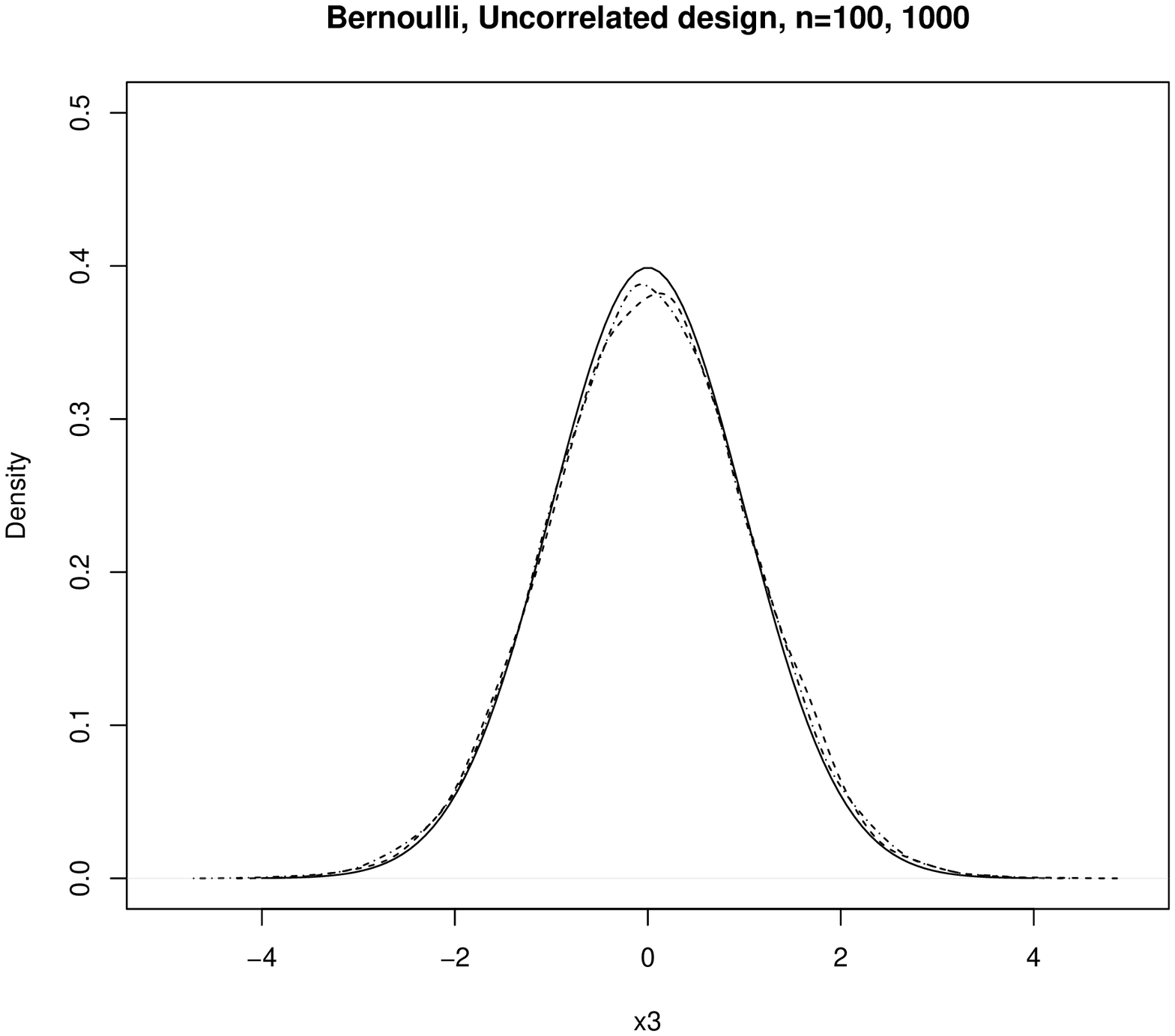}\\
\caption{The density estimate of $U_i$, $[U_i\ U_j]^T$ and the density of $N(0,1)$ and $N(\vec{0},I_2)$ with the Bernoulli model and an uncorrelated design. 
For $i=1,2,3$, the $(i,i)$ panels are the density estimates of $U_i$ for $n=100$(dot-dashed) and $n=1000$(dashed), and the density of $N(0,1)$(solid). 
The (2,1), (3,1) and (3,2) panels are the density estimates of $[U_1\ U_2]^T$, $[U_1\ U_3]^T$ and $[U_2\ U_3]^T$(dashed) for $n=100$ and the density of $N(\vec{0},I_2)$(solid).
The (1,2), (1,3) and (2,3) panels are the density estimates of $[U_1\ U_2]^T$, $[U_1\ U_3]^T$ and $[U_2\ U_3]^T$(dashed) for $n=1000$ and the density of $N(\vec{0},I_2)$(solid).
In each panel, the contour lines of $N(\vec{0},I_2)$ are the same as that of the density estimate. \label{log.simu}}
\end{center}
\end{figure}

\begin{figure}
\begin{center}
\includegraphics[width=45mm,height=45mm]{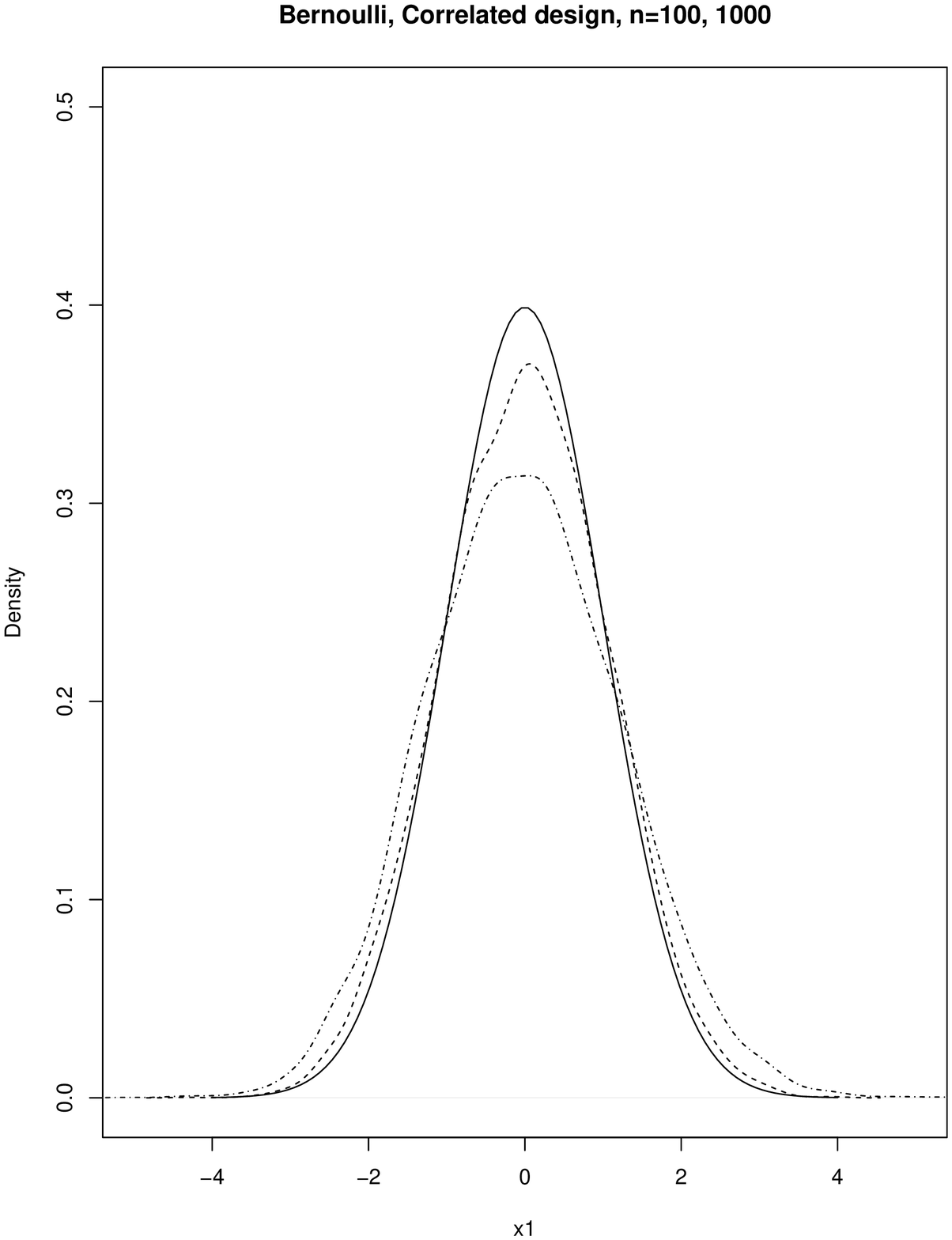}
\includegraphics[width=45mm,height=45mm]{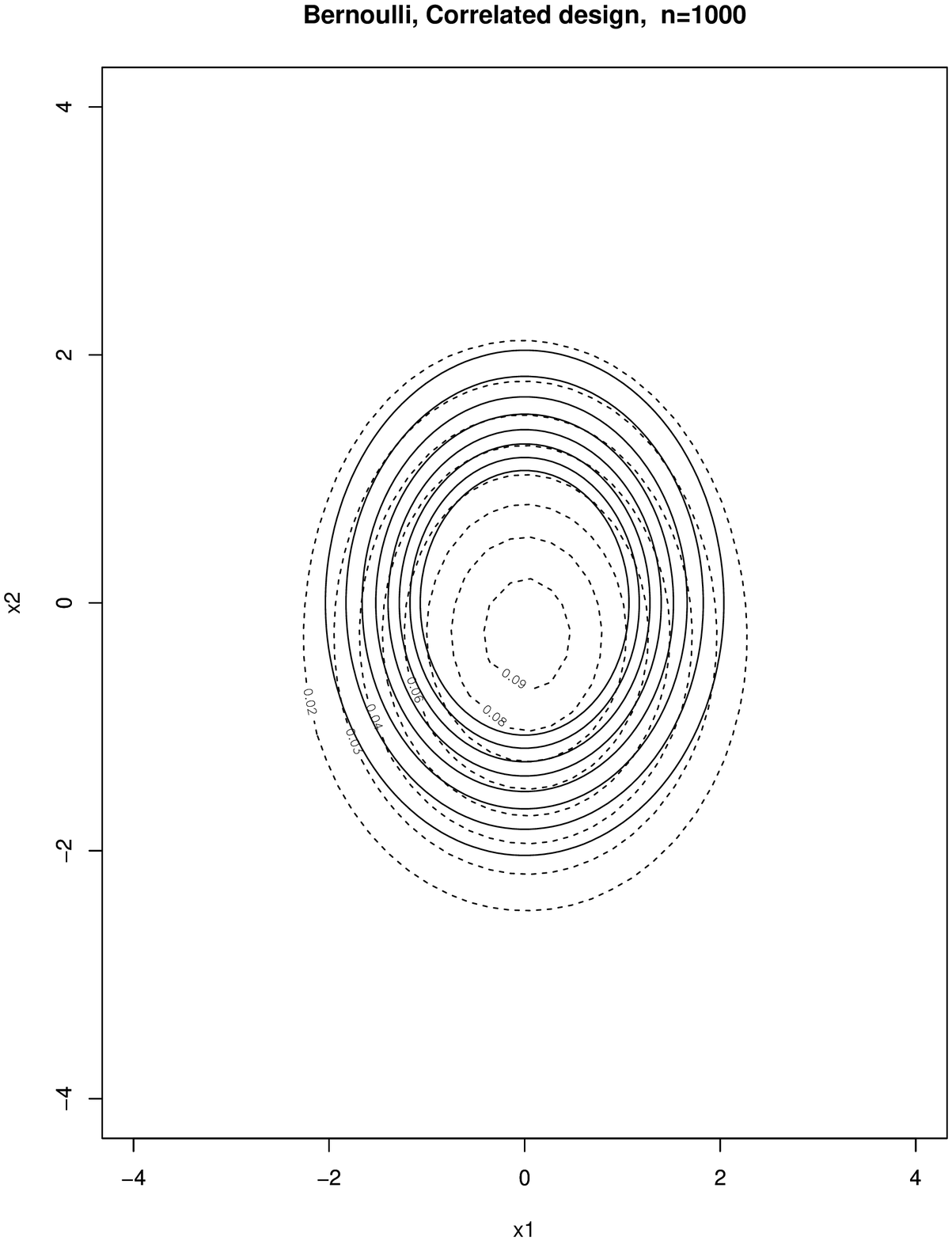}
\includegraphics[width=45mm,height=45mm]{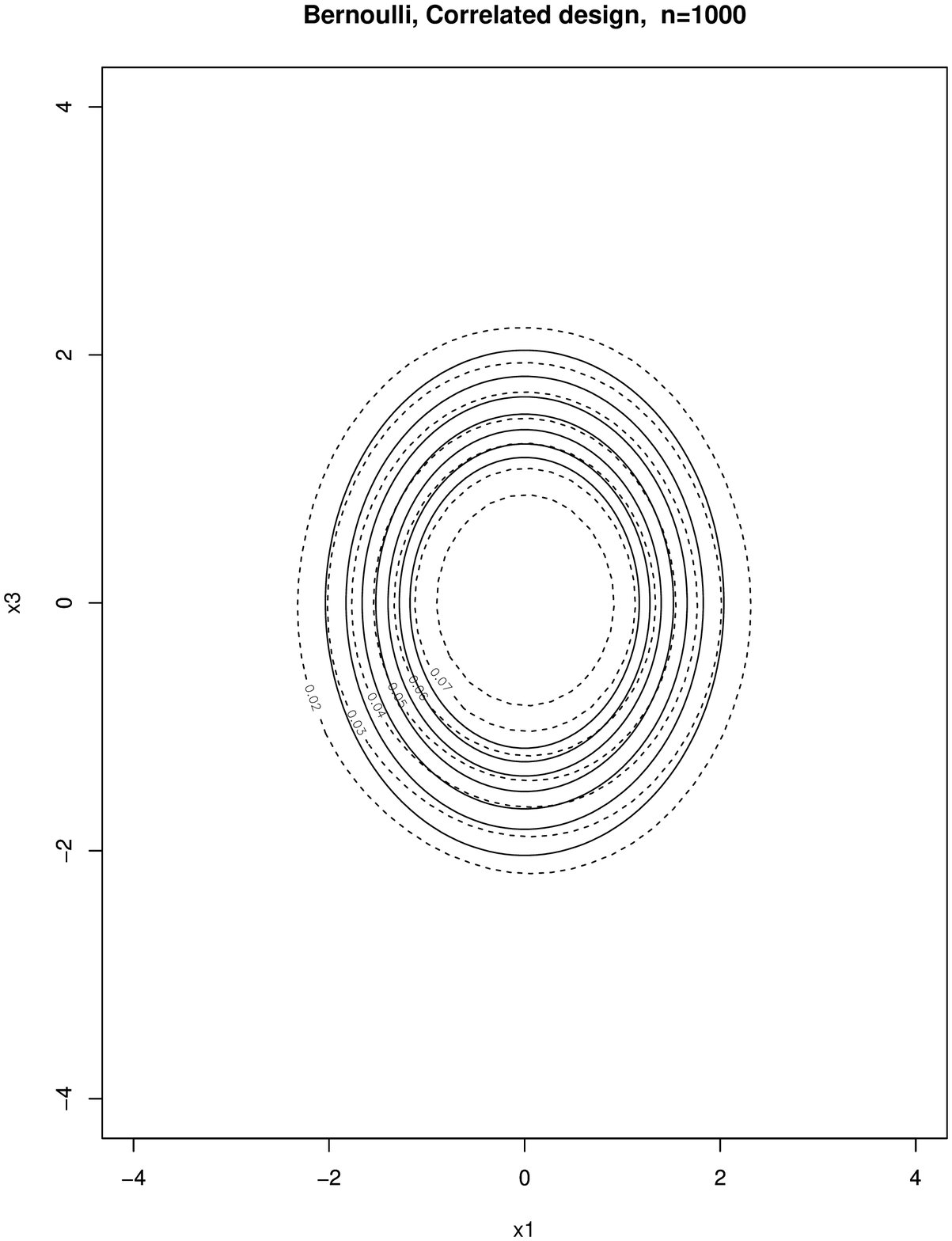}\\
\includegraphics[width=45mm,height=45mm]{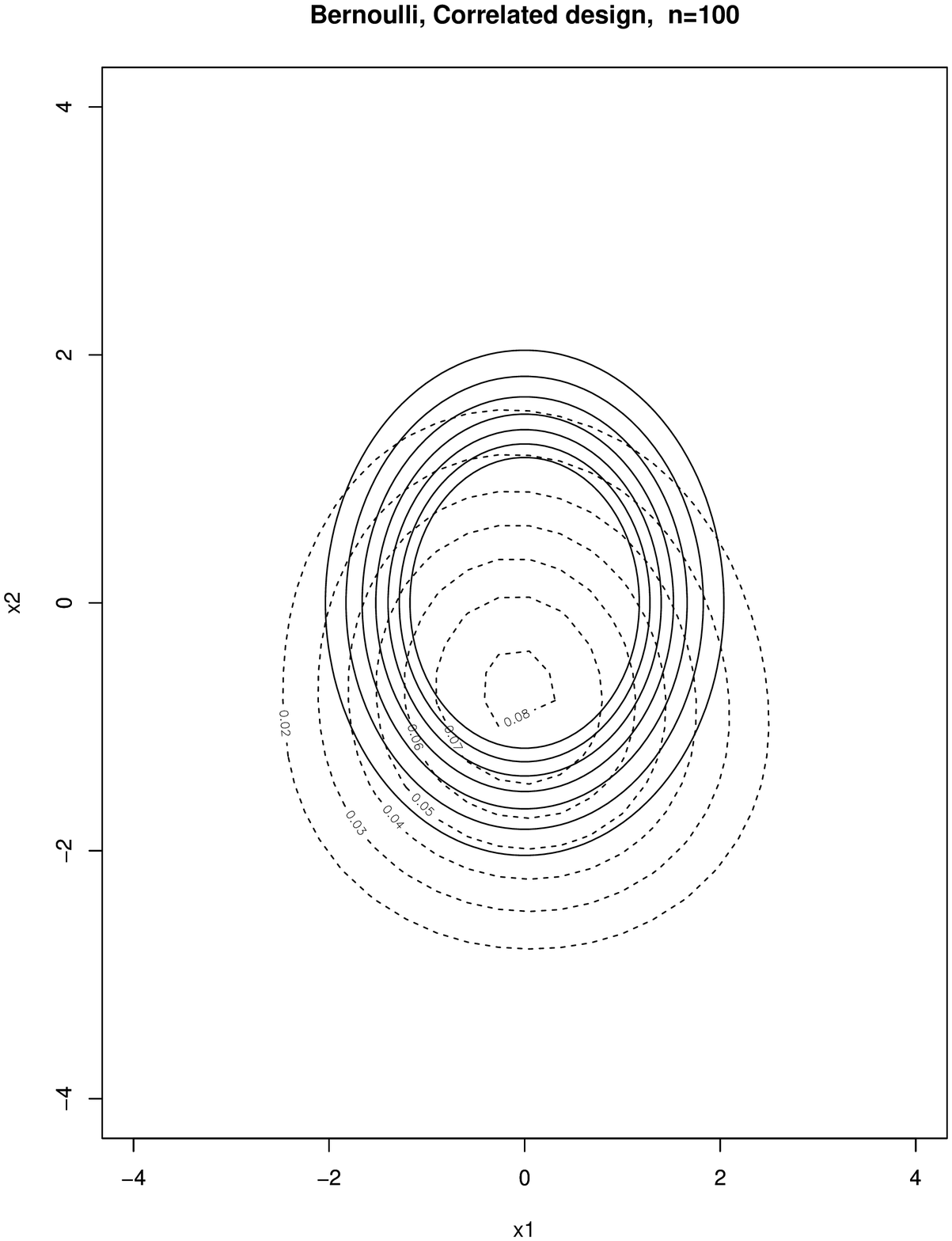}
\includegraphics[width=45mm,height=45mm]{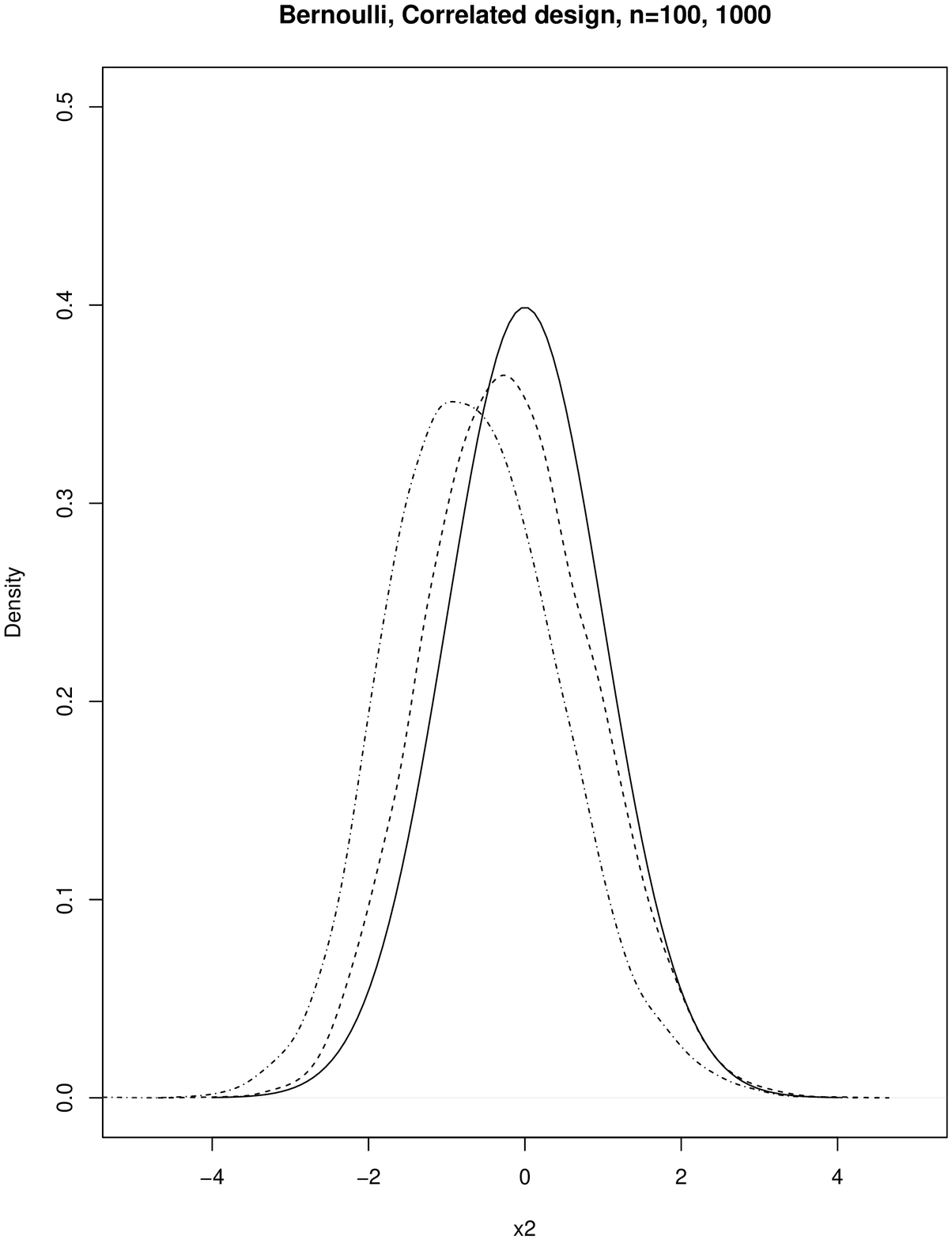}
\includegraphics[width=45mm,height=45mm]{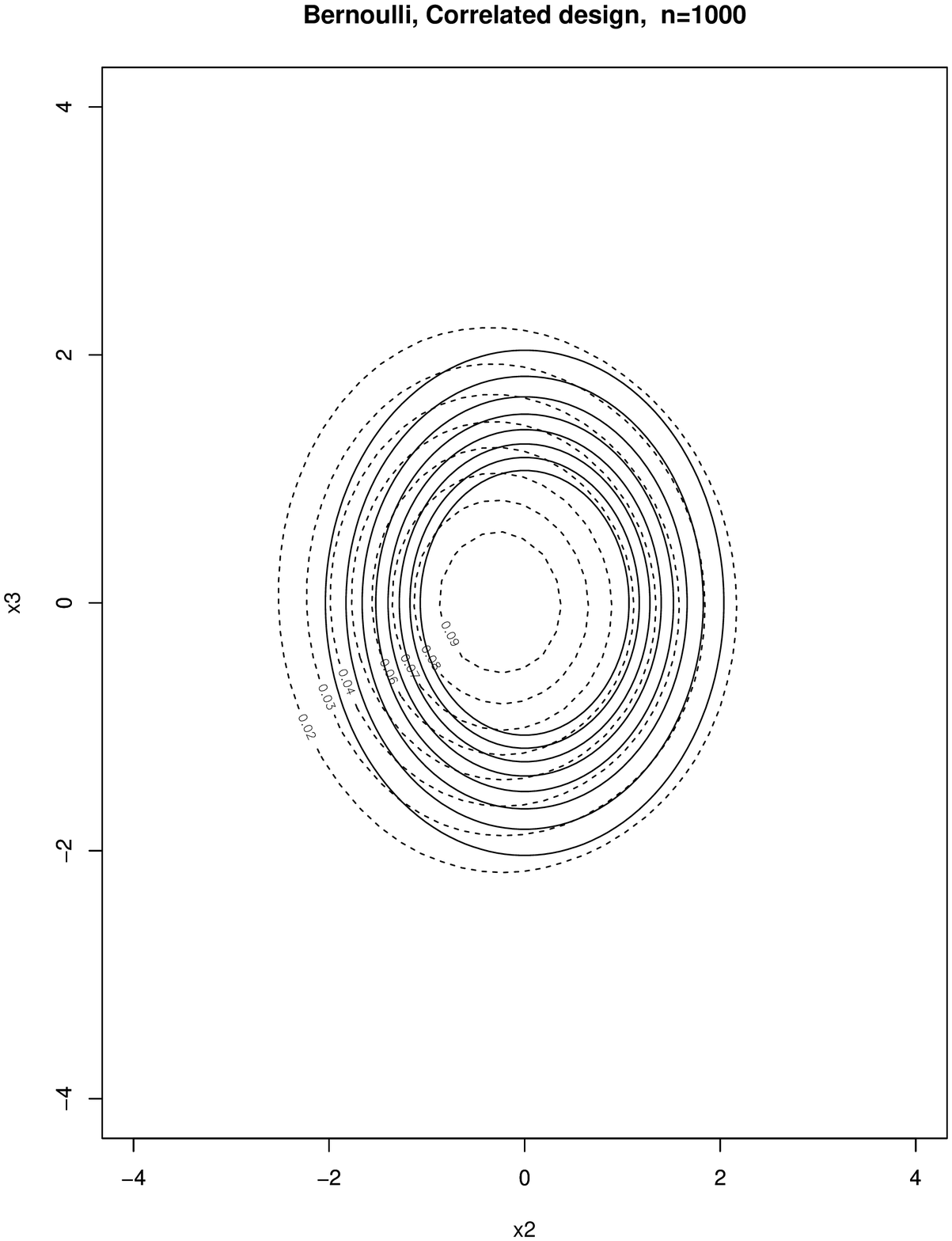}\\
\includegraphics[width=45mm,height=45mm]{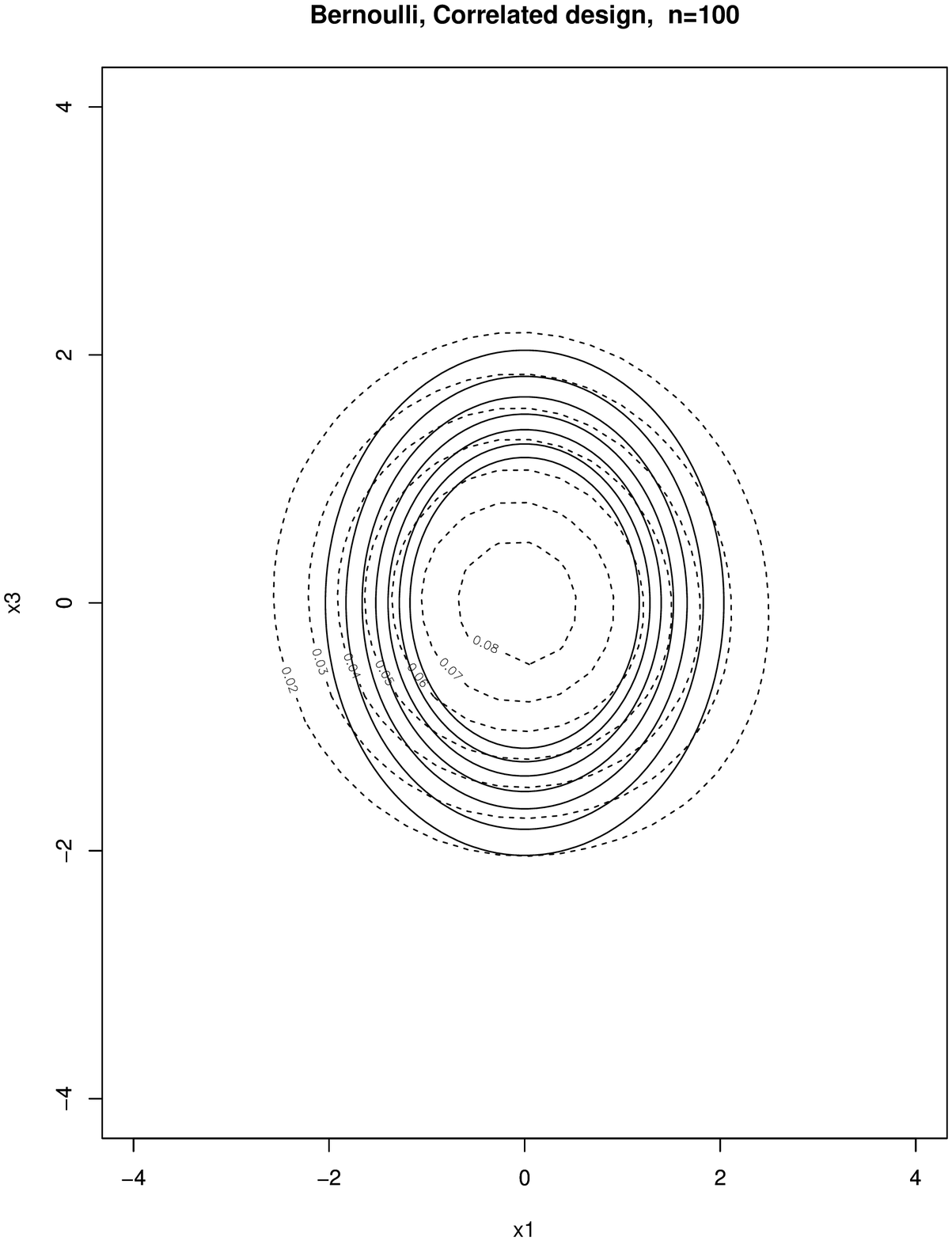}
\includegraphics[width=45mm,height=45mm]{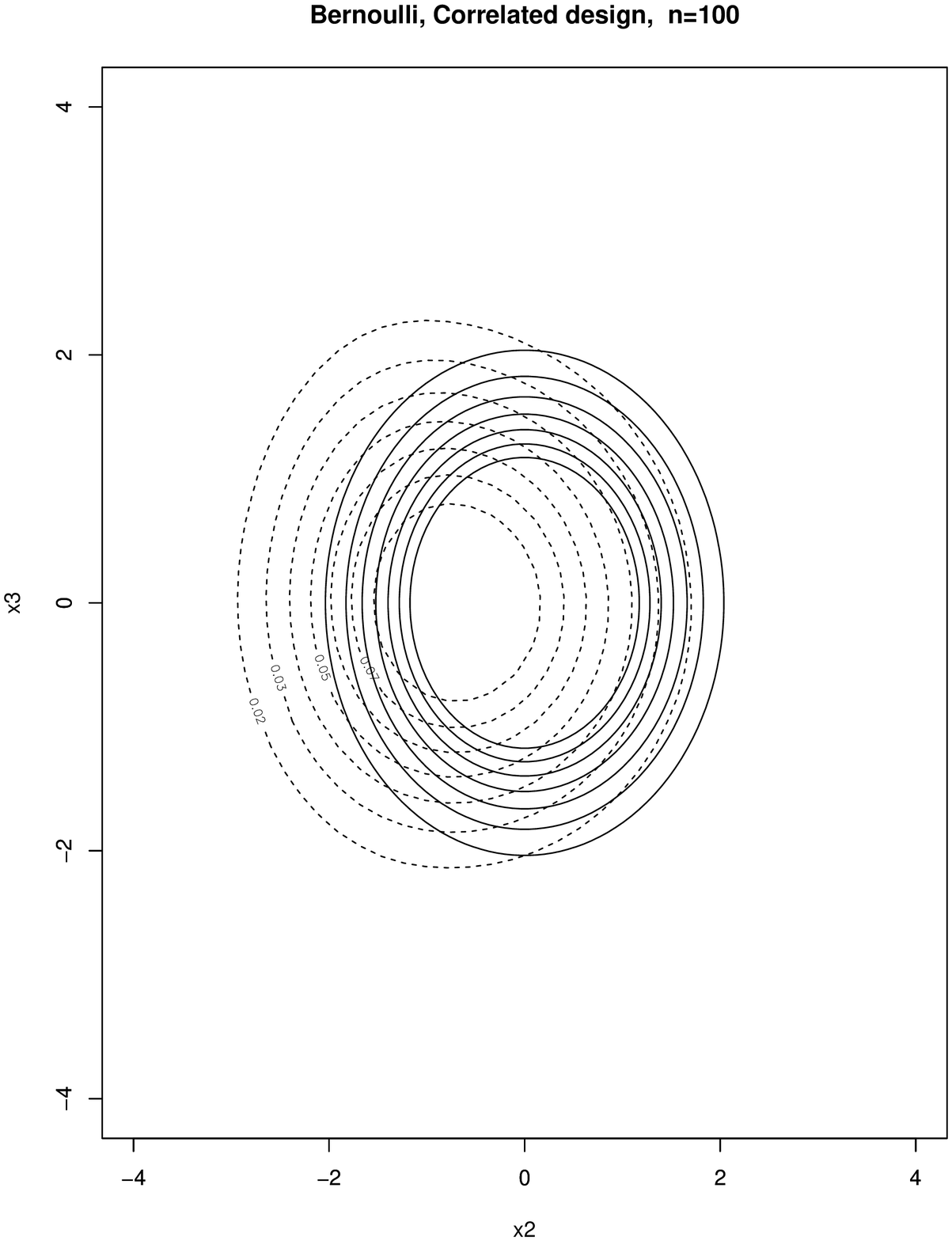}
\includegraphics[width=45mm,height=45mm]{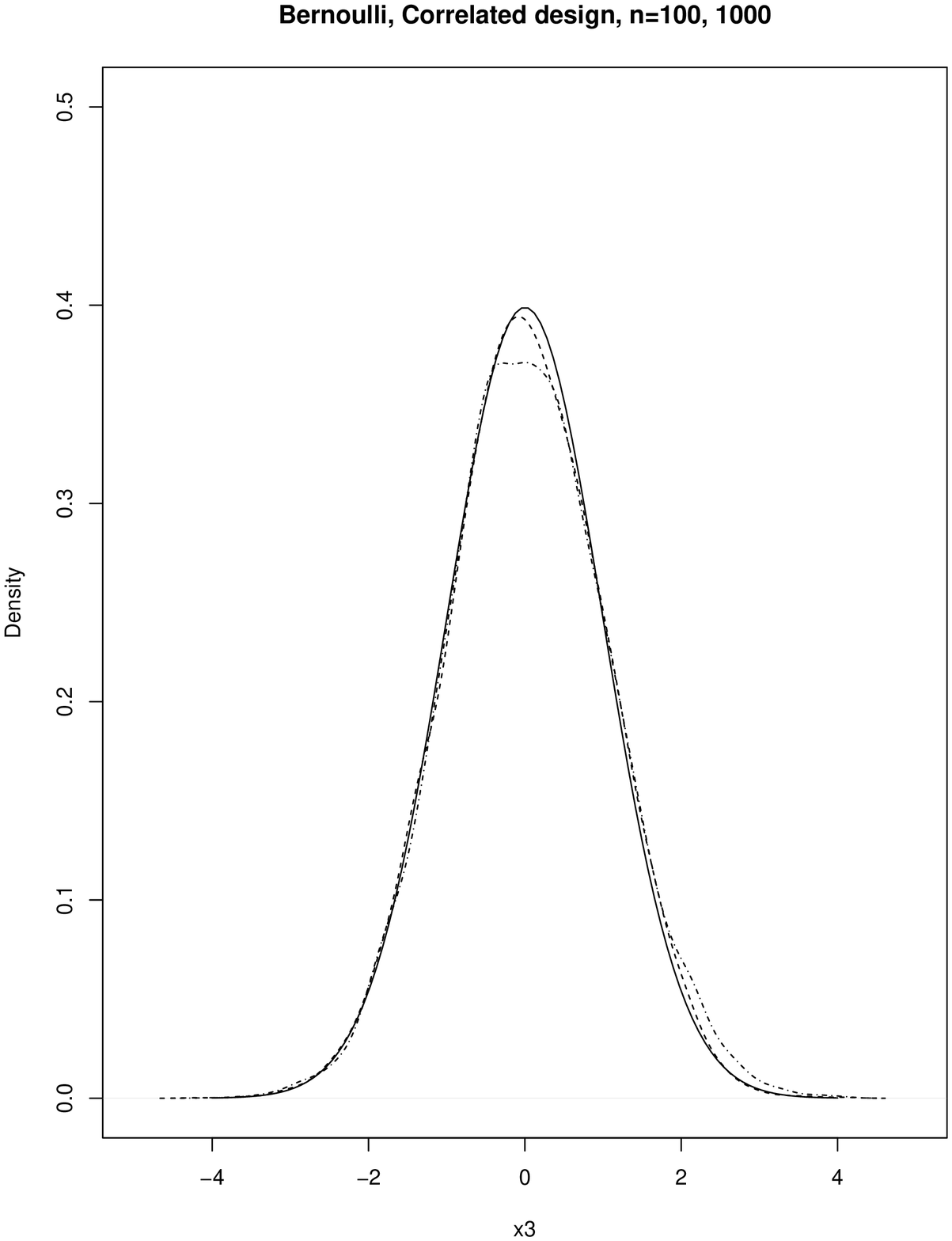}\\
\caption{The density estimate of $U_i$, $[U_i\ U_j]$ and the density of $N(0,1)$ and $N(\vec{0},I_2)$ with the Bernoulli model and the correlated design. 
The description of each panel is the same as in Fig. \ref{log.simu}.  \label{log.simu.cor}}
\end{center}
\end{figure}

In Fig. \ref{log.simu}, the density estimate of (\ref{2norm}) with (i), and the densities of the normal distribution are shown. 
As the sample size increases, the asymptotic normality of the RCPS in Theorem \ref{clt} can be observed numerically. 
We see from (1,1), (2,2) and (3,3) panels that the density estimate becomes close to 0 when $n=1000$. 
When $n=100$, a large correlation between $U_i$ and $U_j$ can be observed. 
However, as $n$ increases, the correlation becomes small.
The results with the correlated design (ii) are drawn in Fig. \ref{log.simu.cor}. 
The density estimate of $U_2$ appears to be far from $N(0,1)$, even when $n=1000$. 
However, we also find that $[U_i\ U_j]^T$ tends to become close to $N_2(\vec{0},I_2)$ as $n\rightarrow \infty$.
We have confirmed that the density estimate with $Y_i\sim {\rm Poisson(\eta(\vec{x}_i))}$ tends to become close to the normal distribution as $n$ increases, though this is not shown in this paper. 
However, the speed of convergence of the density estimate with the Poisson model was somewhat slower than with the Bernoulli model.

\section{Discussion}

This paper showed the asymptotic normality of the penalized spline estimator in the GAM. 
The results of this paper generalize Theorem 1 of Kauermann et al. (2009) and Theorem 2 of Yoshida and Naito (2012). 
The main tools used to prove our Theorems were the spline approximation theories and the asymptotic properties of the band matrices. 
By applying their properties, the asymptotics for penalized splines in other models can be investigated for further study.

In spline smoothing, the determination of smoothing parameters is very important. 
Many researchers have addressed this problem by using grid search methods, such as Mallow's $C_p$, cross-validation and generalized cross-validation. 
Since the computation time of a grid search is dramatically increased when $D>1$, more direct methods would be a useful area of study. 
It may be possible to discuss the selection of smoothing parameters based on the asymptotic properties in this paper.

In recent years, the so-called high-dimensional additive models characterized by $^^ ^^ n<D$" have been studied by many authors such as Meier et al. (2009), Huang et al. (2010) and Fan et al. (2011). 
These previous works are based on unpenalized $B$-spline estimators.
Although it is beyond the scope of this paper, the asymptotics for penalized splines in high dimensional additive models would be interesting to explore.

\section*{Appendix}

For a matrix $X_n=(X_{ij,n})_{ij}$, if $\displaystyle\max_{i,j}\{n^{\alpha}|X_{ij,n}|\}=O_{P}(1)(o_{P}(1))$, then it is written as $X_n=O_{P}(n^{-\alpha}\vec{1}\vec{1}^T)(o_{P}(n^{-\alpha}\vec{1}\vec{1}^T))$. 
When $X_n$ is vector, it is written as $X_n=O_{P}(n^{-\alpha}\vec{1})$. 
We define $W_0=\diag[c^{\prime\prime}(Z\vec{b}_0)]$, 
$G_{j,n}=n^{-1}Z_j^T W_0Z_j$, $G_{i,j,n}=n^{-1}Z_i^T W_0Z_j$ and  $G_{j,i,n}=G_{i,j,n}^T$. 
In the sequel, we use $H_{j,n}=G_{jn}+(\lambda_{jn}/n)\Delta_m^\prime \Delta_m+(\gamma_n/n) I (i,j=1,\cdots,D,i\not=j)$.

We need 3 additional Lemmas as follows.

\begin{lemma}\label{G1}
$G_{j,n}$, $G_{i,j,n}$ and $H_{j,n}$ satisfy $G_{j,n}=O_P(K_n^{-1}\vec{1}\vec{1}^T)$, $G_{i,j,n}=O_P(K_n^{-2}\vec{1}\vec{1}^{T})$ and $H_{j,n}^{-1}=O_P(K_n\vec{1}\vec{1}^T)$. 
Let $A=(a_{ij})_{ij}$ be $(K_n+p)\times (K_n+p)$ matrix. 
Assume that as $K_n\rightarrow \infty$, $A=O_P(K_n^\alpha\vec{1}\vec{1}^T)$.
Then, under the Assumptions, $G_{jn}A=O_P(K_n^{\alpha-1}\vec{1}\vec{1}^T)$, $G_{i,j,n}A=O_P(K_n^{\alpha-2}\vec{1}\vec{1}^{T})$ and 
$H_{j,n}^{-1}A=O_P(K_n^{1+\alpha}\vec{1}\vec{1}^T)$.
\end{lemma}

\begin{lemma}\label{G2}
Let $A_{D,n}$ be $\{D(K_n+p)\}\times \{D(K_n+p)\}$ matrix. 
Assume that as $K_n\rightarrow \infty$, $A_{D,n}=O_P(K_n^{\alpha}\vec{1}\vec{1}^T)$. 
Then, under the Assumptions,  $A_{D,n}H(\vec{b}_0,\lambda_n,\gamma_n)^{-1}=O_P(K_n^{\alpha+1}\vec{1}\vec{1}^T)$
\end{lemma}

\begin{lemma}\label{G3}
Under the assumption, for $j=1,\cdots,D$, $\Delta_m\vec{b}_{j0}=O(K_n^{-m}\vec{1})$.
\end{lemma}

Lemma \ref{G1} can be proven by the properties of the integral of $B$-spline basis and the inverse of band matrices detailed in Claeskens et al. (2009) and Yoshida and Naito (2012). 
Then, Assumption 5 of this paper indicates that the case $K_q<1$ of Claeskens et al. (2009).
The proof of Lemma \ref{G2} is addressed in Yoshida and Naito (2012) by induction for $D$.
Lemma \ref{G3} can be shown from the derivative property of $B$-spline model: $s^{(m)}(x)=\vec{B}^{[p-m]}(x)^T K_n^m\Delta_m\vec{b}$. 
The above equality and Proposition 1 yield $\vec{B}^{[p-m]}(x)^T K_n^m\Delta_m\vec{b}_{j0}=\eta_j^{(m)}(x)(1+o(1))$. 
Since the asymptotic order of $\vec{B}^{[p-m]}(x)^T K_n^m\Delta_m\vec{b}_{j0}$ and each component of $K_n^m\Delta_m\vec{b}_{j0}$ are the same as $O(1)$, Lemma \ref{G3} holds. 
The details are clarified in Section 2 of Claeskens et al. (2009).

\begin{lemma}\label{ex}
Under the Assumptions,
\begin{eqnarray*}
\hat{\vec{b}}-\vec{b}_0=-H(\vec{b}_0,\lambda_n,\gamma_n)^{-1}G(\vec{b}_0,\lambda_n,\gamma_n)+o_P\left(\left\{\left(\frac{\lambda_nK_n^{1-m}}{n}\right)^2+\frac{K_n}{n}\right\}\vec{1}\right).
\end{eqnarray*}
\end{lemma}

\vspace{5mm}

\subsection*{proof of Lemma \ref{ex}}

We use the Taylor expansion of $G(\hat{\vec{b}},\lambda_n,\gamma_n)$ around $\vec{b}_0$, giving
\begin{eqnarray*}
0&=&G(\hat{\vec{b}},\lambda_n,\gamma_n)\\
&=&G(\vec{b}_0,\lambda_n,\gamma_n)+H(\vec{b}_0,\lambda_n,\gamma_n)(\hat{\vec{b}}-\vec{b}_0)\\
&&+ \{H(\vec{b}_0+\Omega(\hat{\vec{b}}-\vec{b}_0),\lambda_n,\gamma_n)-H(\vec{b}_0,\lambda_n,\gamma_n)\}(\hat{\vec{b}}-\vec{b}_0),
\end{eqnarray*}
where
$\Omega={\rm diag}[\omega_{1}\ \cdots\ \omega_{D(K_n+p)}]$ and $\omega_i\in (0,1)$. 
Therefore, $\hat{\vec{b}}-\vec{b}_0$ can be written as 
\begin{eqnarray}
\hat{\vec{b}}-\vec{b}_0
&=&
\left\{-H(\vec{b}_0,\lambda_n,\gamma_n)\right\}^{-1}G(\vec{b}_0,\lambda_n,\gamma_n)\nonumber\\
&&-
H(\vec{b}_0,\lambda_n,\gamma_n)^{-1}\{H(\vec{b}_0+\Omega(\hat{\vec{b}}-\vec{b}_0),\lambda_n,\gamma_n)-H(\vec{b}_0,\lambda_n,\gamma_n)\}(\hat{\vec{b}}-\vec{b}_0)\nonumber\\
&=&
\left\{-H(\vec{b}_0,\lambda_n,\gamma_n)\right\}^{-1}G(\vec{b}_0,\lambda_n,\gamma_n)\nonumber\\
&&+
H(\vec{b}_0,\lambda_n,\gamma_n)^{-1}\left(\frac{1}{n}Z^T\left\{W(\vec{b}_0+\Omega(\hat{\vec{b}}-\vec{b}_0))-W(\vec{b}_0) \right\}Z\right)(\hat{\vec{b}}-\vec{b}_0)
,\label{ex21}
\end{eqnarray}
where $W(\vec{b})=\diag[c^{\prime\prime}(Z(\vec{x}_i)^T\vec{b})]$. 
Furthermore, for $i=1,\cdots,n$, the Taylor expansion yields
\begin{eqnarray*}
c^{\prime\prime}(Z(\vec{x}_i)^T\{\vec{b}_0+\Omega(\hat{\vec{b}}-\vec{b}_0)\})
&=&c^{\prime\prime}(Z(\vec{x}_i)^T\vec{b}_0)\\
&&+
c^{(3)}[Z(\vec{x}_i)^T\vec{b}_0+\theta_iZ(\vec{x}_i)^T\Omega(\hat{\vec{b}}-\vec{b}_0)]Z(\vec{x}_i)^T\Omega(\hat{\vec{b}}-\vec{b}_0),
\end{eqnarray*}
where $\theta_i\in(0,1)$. 
Hence from Proposition 1, we obtain 
\begin{eqnarray*}
W(\vec{b}_0+\Omega(\hat{\vec{b}}-\vec{b}_0))-W(\vec{b}_0)
&=&\diag[c^{(3)}(Z(\vec{x}_i)^T\vec{b}_0+\theta_iZ(\vec{x}_i)^T\Omega(\hat{\vec{b}}-\vec{b}_0))Z(\vec{x}_i)^T\Omega(\hat{\vec{b}}-\vec{b}_0)]\\
&=&\diag[c^{(3)}(\eta(\vec{x}_i))Z(\vec{x}_i)^T\Omega(\hat{\vec{b}}-\vec{b}_0)](1+o(1))\\
&\equiv& R(\hat{\vec{b}}).
\end{eqnarray*}
 
For simplicity, we rewrite $G=G(\vec{b}_0,\lambda_n,\gamma_n)$ and $H=-H(\vec{b}_0,\lambda_n,\gamma_n)$. 
Then (\ref{ex21}) can be written as 
\begin{eqnarray}
\hat{\vec{b}}-\vec{b}_0=H^{-1}G+H^{-1}\left(\frac{1}{n}Z^T R(\hat{\vec{b}})Z\right)(\hat{\vec{b}}-\vec{b}_0). \label{ex22}
\end{eqnarray}
We now prove 
\begin{eqnarray}
H^{-1}\left(\frac{1}{n}Z^T R(\hat{\vec{b}})Z\right)(\hat{\vec{b}}-\vec{b}_0)
=
o_P\left(\left\{\left(\frac{\lambda_nK_n^{1-m}}{n}\right)^2+\frac{K_n}{n}\right\}\vec{1}\right). \label{purp1}
\end{eqnarray}
From (\ref{ex22}), the left hand side of (\ref{purp1}) can be evaluated as 
\begin{eqnarray*}
&&H^{-1}\left(\frac{1}{n}Z^T R(\hat{\vec{b}})Z\right)(\hat{\vec{b}}-\vec{b}_0)\\
&&=
H^{-1}\left(\frac{1}{n}Z^T R(\hat{\vec{b}})Z\right)H^{-1}G+\left\{H^{-1}\left(\frac{1}{n}Z^T R(\hat{\vec{b}})Z\right)\right\}^2(\hat{\vec{b}}-\vec{b}_0).
\end{eqnarray*}
First we show the asymptotic order of $R(\hat{\vec{b}})$.
The $i$th component of $R(\hat{\vec{b}})$ can be written by (\ref{ex22}) as  
\begin{eqnarray*}
&&c^{(3)}(\eta(\vec{x}_i))Z(\vec{x}_i)^T\Omega(\hat{\vec{b}}-\vec{b}_0)(1+o(1))\\
&&= c^{(3)}(\eta(\vec{x}_i))Z(\vec{x}_i)^T\Omega\left\{H^{-1}G+H^{-1}\left(\frac{1}{n}Z^T R(\hat{\vec{b}})Z\right)(\hat{\vec{b}}-\vec{b}_0)\right\}(1+o(1)).
\end{eqnarray*}
By calculating the expectation and the square root of variance of each component of $G$, we obtain with Lemma \ref{G3}
$$
G=O_P\left(\left\{\frac{\lambda_n}{nK_n^{m}}+\frac{1}{\sqrt{nK_n}}\right\}\vec{1}\right).
$$
Therefore Lemma \ref{G2} yields that for $\vec{z}\in[0,1]^D$, 
\begin{eqnarray*}
Z(\vec{z})^T\Omega H^{-1}G=o_P\left(\left\{\frac{\lambda_nK_n^{1-m}}{n}+\sqrt{\frac{K_n}{n}}\right\}\right).
\end{eqnarray*}
Since $c^{(3)}(\eta(\vec{x}))$ is bounded near $\eta(\vec{x})$ for $\vec{x}\in[0,1]^D$, we have with tedious but easy calculation that
\begin{eqnarray}
&&\left|c^{(3)}(\eta(\vec{x}_i))Z(\vec{x}_i)^T\Omega\left\{H^{-1}G+H^{-1}\left(\frac{1}{n}Z^T R(\hat{\vec{b}})Z\right)(\hat{\vec{b}}-\vec{b}_0)\right\}\right|(1+o(1))\nonumber\\
&&\leq \sup_{\vec{z}\in[0,1]^D}\left|c^{(3)}(\eta(\vec{z}))Z(\vec{z})^T\Omega\left\{H^{-1}G+H^{-1}\left(\frac{1}{n}Z^T R(\hat{\vec{b}})Z\right)(\hat{\vec{b}}-\vec{b}_0)\right\}\right|(1+o(1))\nonumber\\
&&=o_P\left(\frac{\lambda_nK_n^{1-m}}{n}+\sqrt{\frac{K_n}{n}}\right).\label{Rcal}
\end{eqnarray}
Then Lemmas \ref{G1} and \ref{G2} and (\ref{Rcal}) yield
\begin{eqnarray*}
H^{-1}\left(\frac{1}{n}Z^T R(\hat{\vec{b}})Z\right)H^{-1}G
&=&
H^{-1}\left(\frac{1}{n}Z^T R(\hat{\vec{b}})Z\right)O_P\left(\left\{\frac{\lambda_nK_n^{1-m}}{n}+\sqrt{\frac{K_n}{n}}\right\}\right)\\
&=&
H^{-1}o_P\left(\frac{1}{K_n}\left\{\frac{\lambda_nK_n^{1-m}}{n}+\sqrt{\frac{K_n}{n}}\right\}^2\right)\\
&=&
o_P\left(\left\{\frac{\lambda_nK_n^{1-m}}{n}+\sqrt{\frac{K_n}{n}}\right\}^2\vec{1}\right).
\end{eqnarray*}
Further we get with simple calculation
\begin{eqnarray*}
\left\{H^{-1}\left(\frac{1}{n}Z^T R(\hat{\vec{b}})Z\right)\right\}^2(\hat{\vec{b}}-\vec{b}_0)=o_P\left(\left\{\frac{\lambda_nK_n^{1-m}}{n}+\sqrt{\frac{K_n}{n}}\right\}^2\vec{1}\right).
\end{eqnarray*}
This implies (\ref{purp1}) and completes the proof.  
$\Box$

\vspace{5mm}

\subsection*{proof of Proposition \ref{app}}

Barrow and Smith (1978) showed that for $j=1,\cdots,D$, there exists $\vec{b}_j^*\in\mathbb{R}^{K_n+p}$ such that 
\begin{eqnarray*}
\sup_{z\in(0,1)}\left|\eta_j(z)+b_{j,a}(z)-\vec{B}(z)^T\vec{b}_j^*\right|
=o(K_n^{-(p+1)}).
\end{eqnarray*} 
Let $\eta_j^*(z)=\vec{B}(z)^T\vec{b}_j^*$, $\eta^*(\vec{x})=\sum_{j=1}^D \eta^*_j(x_j)$, $\eta_0(\vec{x})=\sum_{j=1}^D \eta_{j0}(x_j)$ and $b_a(\vec{x})=\sum_{j=1}^D b_{j,a}(x_j)$. 
We now prove that 
\begin{eqnarray}
\vec{b}_{j0}-\vec{b}_j^*=o(K_n^{-(p+1)}\vec{1}),\ \ j=1,\cdots,D. \label{pur}
\end{eqnarray}
Since the asymptotic order of $\eta_{j0}(x_j)-\eta_j^*(x_j)$ and that of $\vec{b}_0-\vec{b}^*$ are the same, if (\ref{pur}) is satisfied, 
we obtain for any $x_j\in(0,1)$,
$|\eta_{j0}(x_j)-\eta_j^*(x_j)|=o(K_n^{-(p+1)})$ hence Proposition 1 holds.

From the definition of $\vec{b}_0$, we have 
\begin{eqnarray}
\frac{1}{n}\sum_{i=1}^n E\left[\left.\log \frac{f(Y_i|\vec{x}_i,\eta)}{f(Y_i|\vec{x}_i,\vec{b}_0)}\right|\vec{X}_n\right]
\leq 
\frac{1}{n}\sum_{i=1}^n E\left[\left.\log \frac{f(Y_i|\vec{x}_i,\eta)}{f(Y_i|\vec{x}_i,\vec{b}^*)}\right|\vec{X}_n\right],\label{ineq}
\end{eqnarray}
where $\vec{b}^*=((\vec{b}^*_1)^T\ \cdots\ (\vec{b}^*_D)^T)^T$.
The Taylor expansion to $c(\eta^*(\vec{x}_i))$ around $\eta(\vec{x}_i)$ yields 
\begin{eqnarray*}
&&\frac{1}{n}\sum_{i=1}^n E\left[\left.\log \frac{f(Y_i|\vec{x}_i,\eta)}{f(Y_i|\vec{x}_i,\vec{b}^*)}\right|\vec{X}_n\right]\\
&&=
\frac{1}{n}\sum_{i=1}^n\left[c^\prime(\eta(\vec{x}_i))\{\eta(\vec{x}_i)-\eta^*(\vec{x}_i)\}-\{c(\eta(\vec{x}_i))-c (\eta^*(\vec{x}_i))\}\right]\nonumber\\
&&=
\frac{1}{2n}\sum_{i=1}^n\left[\{\eta(\vec{x}_i)-\eta^*(\vec{x}_i)\}^2c^{\prime\prime} (\eta(\vec{x}_i))(1+o(1))\right]\nonumber\\
&&=
\frac{1}{2n}\sum_{i=1}^n\left[b_{a}(\vec{x}_i)^2c^{\prime\prime} (\eta(\vec{x}_i)(1+o(1))\right]\\
&&= O(K_n^{-2(p+1)}).
\end{eqnarray*}
Therefore we obtain
$|\eta(\vec{x}_i)-\eta_{0}(\vec{x}_i)|= o(1)$, by which 
\begin{eqnarray}
\frac{1}{n}\sum_{i=1}^n E\left[\left.\log \frac{f(Y_i|\vec{x}_i,\eta)}{f(Y_i|\vec{x}_i,\vec{b}_0)}\right|\vec{X}_n\right]
&=&\frac{1}{n}\sum_{i=1}^n\left[\{\eta(\vec{x}_i)-\eta_{0}(\vec{x}_i)\}^2c^{\prime\prime} (\eta(\vec{x}_i))(1+o(1))\right]\nonumber\\
&=&
\frac{1}{n}(\vec{\eta}-Z\vec{b}_0)^T W (\vec{\eta}-Z\vec{b}_0), \label{b0def}
\end{eqnarray}
where $W=\diag[c^{\prime\prime}(\eta(\vec{x}_i))(1+o(1))]$ and $\vec{\eta}=(\eta(\vec{x}_1)\ \cdots\ \eta(\vec{x}_n))^T$. 
It is easy to show that $\vec{b}_0$ satisfies 
\begin{eqnarray}
\frac{1}{n}Z^T W Z \vec{b}_0=\frac{1}{n}Z^T W \vec{\eta}
\end{eqnarray}
since $\vec{b}_0$ is the minimizer of (\ref{b0def}).
Further, from the definition of $\vec{b}^*$, we have  
\begin{eqnarray*}
\vec{\eta}=Z\vec{b}^* -\vec{B}_a +o(K_n^{-(p+1)}\vec{1}),
\end{eqnarray*}
where $\vec{B}_a=(b_a(\vec{x}_1)\ \cdots\ b_a(\vec{x}_n))^T$. 
Hence, we obtain 
\begin{eqnarray}
\frac{1}{n}Z^T W Z (\vec{b}_0-\vec{b}^*)=-\frac{1}{n}Z^T W \{\vec{B}_a+o(K_n^{-(p+1)}\vec{1})\}. \label{dist2}
\end{eqnarray}
By noting $Z=[Z_1\ \cdots\ Z_D]$, the $k$th component of first $(K_n+p)$ block of $n^{-1}Z^T W\vec{B}_a$ can be calculated as 
\begin{eqnarray}
\left(n^{-1}Z_1^T W\vec{B}_a\right)_k
&=&
\frac{1}{n}\sum_{i=1}^n c^{\prime\prime}(\eta(\vec{x}_i))B_{-p+k}(x_{i1})b_a(\vec{x}_{i})(1+o(1))\nonumber\\
&=&
\sum_{j=1}^D \left[\frac{1}{n}\sum_{i=1}^n c^{\prime\prime}(\eta(\vec{x}_i))B_{-p+k}(x_{i1})b_{j,a}(x_{ij})\right](1+o(1))\nonumber\\
&=&
\sum_{j=1}^D \int_{[0,1]^D} c^{\prime\prime}(\eta(\vec{x}))B_{-p+k}(x_{1})b_{j,a}
(x_{j})dP(\vec{x})(1+o(1))\nonumber\\
&=&o(K_n^{-(p+2)}).\label{brsp}
\end{eqnarray}
Here the last equality in (\ref{brsp}) can be obtained by mimicking the proof of Lemma 10 of Agarwal and Studden (1980). 
Similarly since the row sum of $n^{-1}Z^T W$ has an order $O(K_n^{-1})$, we get (\ref{dist2}) as 
\begin{eqnarray*}
\frac{1}{n}Z^T W Z (\vec{b}_0-\vec{b}^*)=-\frac{1}{n}Z^T W \{\vec{B}_a+o(K_n^{-(p+1)}\vec{1})\}=o(K_n^{-(p+2)}\vec{1}).
\end{eqnarray*}
From Lemma \ref{G1}, we have $n^{-1}Z^T W Z=O(K_n^{-1}\vec{1}\vec{1}^T)$. 
Therefore 
$$
\vec{b}_0-\vec{b}^*=o(K_n^{-(p+1)}\vec{1})
$$
and (\ref{pur}) can be proven. 
$\Box$

\vspace{5mm}

To complete the proof of Theorem 1, first we will obtain (A) the asymptotic form of $\hat{\vec{b}}_D-\vec{b}_{D0}$. And then (B) we will derive the asymptotic form of $E[\hat{\vec{b}}_D-\vec{b}_{D0}|\vec{X}_n]$ and $V[\hat{\vec{b}}_D|\vec{X}_n]$.
Similar argument will be applied to $\hat{\vec{b}}_j-\vec{b}_{j0} (j=1,\cdots,D-1).$

\vspace{5mm}

\subsection*{proof of Theorem \ref{mv}}

First we aim to show (A).
From Lemma \ref{ex}, we obtain 
\begin{eqnarray*}
\hat{\vec{b}}-\vec{b}_0
&=&
\left\{-H(\vec{b}_0,\lambda_n,\gamma_n)\right\}^{-1}G(\vec{b}_0,\lambda_n,\gamma_n)+R_n(\hat{\vec{b}})\\
&=&
\left(\frac{1}{n}Z^T W_0 Z+ \frac{1}{n}Q_m(\lambda_n)+\frac{\gamma_n}{n} I\right)^{-1}G(\vec{b}_0,\lambda_n,\gamma_n)+ R_n(\hat{\vec{b}}), 
\end{eqnarray*}
where 
$$
R_n(\hat{\vec{b}})=o_P\left(\left\{\left(\frac{\lambda_nK_n^{1-m}}{n}\right)^2+\frac{K_n}{n}\right\}\vec{1}\right).
$$
Let $M_D=n^{-1}(Z^T W_0 Z+Q_m(\lambda_n)+\gamma_n I)$ and let $\Lambda_{j,\gamma}=n^{-1}(Z_j^T W_{0} Z_j+\lambda_{jn}\Delta_m^\prime \Delta_m +\gamma_n I)$.  
Then, $M_{D}$ can be written by using $M_{D-1}$ as 
\begin{eqnarray}
M_{D}
=
\left[
\begin{array}{cccc}
\Lambda_{1,\gamma}&G_{1,2,n}&\cdots&G_{1,D,n}\\
G_{2,1,n}&\Lambda_{2,\gamma}&\cdots&G_{2,D,n}\\
\vdots&\vdots &\ddots&\vdots\\
G_{D,1,n}&\cdots&\cdots&\Lambda_{D,\gamma}
\end{array}
\right]
=
\left[
\begin{array}{cc}
M_{D-1}&R^T\\
R&\Lambda_{D,\gamma}
\end{array}
\right],\label{Mmat1}
\end{eqnarray}
where $R=[G_{D,1,n}\ \cdots\ G_{D,D-1,n}]$.
From the result of partitioned matrix (see, Horn and Johnson (1985)), we have
\begin{eqnarray*}
M_D^{-1}=\left[
\begin{array}{cc}
M^{-1}_{D-1}+M^{-1}_{D-1}R^T V^{-1}RM^{-1}_{D-1}&-M^{-1}_{D-1}R^T V^{-1}\\
-V^{-1}RM^{-1}_{D-1}&V^{-1}
\end{array}
\right],
\end{eqnarray*}
where $V=\Lambda_{D,\gamma}-RM^{-1}_{D-1}R^T$.
Let $G_{(-D)}(\vec{b}_0,\lambda_n,\gamma_n)$ and $G_D(\vec{b}_0,\lambda_n,\gamma_n)$ be the first $(D-1)(K_n+p)$th subvector and last $(K_n+p)$th subvector of $G(\vec{b}_0,\lambda_n,\gamma_n)$.
Then, 
\begin{eqnarray*}
\hat{\vec{b}}-\vec{b}_0&\equiv&\left[
\begin{array}{cc}
\hat{\vec{b}}_{(-D)}-\vec{b}_{(-D)0}\\
\hat{\vec{b}}_D-\vec{b}_{D0}
\end{array}
\right]\\
&=&
\left[
\begin{array}{cc}
M^{-1}_{D-1}+M^{-1}_{D-1}R^T V^{-1}RM^{-1}_{D-1}&-M^{-1}_{D-1}R^T V^{-1}\\
-V^{-1}RM^{-1}_{D-1}&V^{-1}
\end{array}
\right]
\left[
\begin{array}{c}
G_{(-D)}(\vec{b}_0,\lambda_n,\gamma_n)\\
G_D(\vec{b}_0,\lambda_n,\gamma_n)
\end{array}
\right]\\
&&\quad + R_n(\hat{\vec{b}}),
\end{eqnarray*}
from which we have 
\begin{eqnarray*}
\hat{\vec{b}}_D-\vec{b}_{0,D}
&=&
V^{-1}G_D(\vec{b}_0,\lambda_n,\gamma_n)-V^{-1}RM^{-1}_{D-1}G_{(-D)}(\vec{b}_0,\lambda_n,\gamma_n)+R_n\\
&=&
\Lambda_{D,\gamma}^{-1}G_D(\vec{b}_0,\lambda_n,\gamma_n) + v_n(\vec{b}_0)+R_{D,n}(\hat{\vec{b}}),
\end{eqnarray*}
where 
\begin{eqnarray*}
v_n(\vec{b}_0)
&=&
-V^{-1}RM^{-1}_{D-1}G_{(-D)}(\vec{b}_0,\lambda_n,\gamma_n)+\left\{V^{-1}-\Lambda_{D,\gamma}^{-1}\right\}G_D(\vec{b}_0,\lambda_n,\gamma_n)
\end{eqnarray*}
and $R_{D,n}(\hat{\vec{b}})$ is last $(K_n+p)$th subvector of $R_n(\hat{\vec{b}})$. 

In following, we shall start to show (B).   
The expectation of $\hat{\vec{b}}_D-\vec{b}_{0,D}$ can be written as 
\begin{eqnarray*}
E[\hat{\vec{b}}_D-\vec{b}_{0,D}|\vec{X}_n]
&=&
\Lambda_{D,\gamma}^{-1}E[G_D(\vec{b}_0,\lambda_n,\gamma_n)|\vec{X}_n] + E[v_n(\vec{b}_0)|\vec{X}_n]+E[R_{D,n}(\hat{\vec{b}})|\vec{X}_n].
\end{eqnarray*}
First, $E[R_{D,n}(\hat{\vec{b}})|\vec{X}_n]=o_P(\lambda_nK_n^{1-m}n^{-1})$ is satisfied.
In the sequel, because 
\begin{eqnarray*}
E[G(\vec{b}_0,0,0)|\vec{X}_n]
&=&
\frac{1}{n}\sum_{i=1}^n E\left[\left.\frac{\partial}{\partial \vec{b}} \log f(Y_i|x_i,\vec{b}_0)\right|\vec{X}_n\right]=0,
\end{eqnarray*}
we have with Lemma \ref{G3}
\begin{eqnarray*}
E[G(\vec{b}_0,\lambda_n,\gamma_n)|\vec{X}_n]&=&E[G(\vec{b}_0,0,0)|\vec{X}_n]-\frac{1}{n}Q_m(\lambda_n)\vec{b}_0-\frac{\gamma_n}{n} \vec{b}_0\\
&=&
O_P\left(\frac{\lambda_nK_n^{-m}}{n}\vec{1}\right).
\end{eqnarray*}
From Lemmas \ref{G1} and \ref{G2}, on the other hand, we obtain 
\begin{eqnarray*}
V^{-1}-\Lambda_{D,\gamma,n}^{-1}
&=&
\Lambda_{D,\gamma}^{-1}(I-RM^{-1}_{D-1}R^T \Lambda_{D,\gamma}^{-1})^{-1}-\Lambda_{D,\gamma}^{-1}\\
&=&
\Lambda_{D,\gamma}^{-1}RM^{-1}_{D-1}R^T \Lambda_{D,\gamma}^{-1}(I-RM^{-1}_{D-1}R^T \Lambda_{D,\gamma}^{-1})^{-1}\\
&=&
O_P\left(K_n\frac{1}{K_n^2}K_n\frac{1}{K_n^2}K_n\vec{1}\vec{1}^T\right)\\
&=&
O_P\left(K_n^{-1}\vec{1}\vec{1}^T\right).
\end{eqnarray*}
Therefore we have with straightforward calculation
\begin{eqnarray*}
&&E[v_n(\vec{b}_0)|\vec{X}_n]\\
&&=
-V^{-1}RM^{-1}_{D-1}E[G_{(-D)}(\vec{b}_0,\lambda_n,\gamma_n)]+\left\{V^{-1}-\Lambda_{D,\gamma,n}^{-1}\right\}E[G_D(\vec{b}_0,\lambda_n,\gamma_n)]\\
&&=
O_P\left(K_n\frac{1}{K_n^2}K_n\frac{\lambda_nK_n^{-m}}{n}\vec{1}\right)
+
O_P\left(\frac{\lambda_nK_n^{-(m+1)}}{n}\vec{1}\right)\\
&&=
O_P\left(\frac{\lambda_nK_n^{-m}}{n}\vec{1}\right).
\end{eqnarray*}

Above calculations are combined into  
\begin{eqnarray*}
E[\hat{\vec{b}}_D-\vec{b}_{D0}|\vec{X}_n]
&=&
-\frac{\lambda_{Dn}}{n}\Lambda_{D,\gamma}^{-1}\Delta_m^\prime \Delta_m\vec{b}_{D0}-\frac{\gamma_n}{n}\Lambda_{D,\gamma,n}^{-1}\vec{b}_{D0}+o_P\left(\frac{\lambda_nK_n^{1-m}}{n}\vec{1}\right)\\
&=&
-\frac{\lambda_{Dn}}{n}\Gamma_D(\lambda_{Dn})^{-1}\Delta_m^\prime \Delta_m\vec{b}_{D0}+o_P\left(\frac{\lambda_nK_n^{1-m}}{n}\vec{1}\right).
\end{eqnarray*}
Here we have used the fact $\Lambda_{D,\gamma}^{-1}=\Gamma_D(\lambda_{Dn})^{-1}(I+o_P(\vec{1}\vec{1}^T))$ and $(\gamma_n/n)\Gamma_D(\lambda_{Dn})^{-1}\vec{b}_{D0}=o(\lambda_nK_n^{1-m}n^{-1}\vec{1})$.
Hence, we finally obtain 
\begin{eqnarray*}
E[\hat{\eta}_{D,\gamma}(x_D)-\eta_{D,0}(x)|\vec{X}_n]
&=&
E[\vec{B}(x_D)^T (\hat{\vec{b}}_D-\vec{b}_{D0})|\vec{X}_n]\\
&=&
b_{D,\lambda}(x_D)+o_P(\lambda_nK_n^{1-m}n^{-1}).
\end{eqnarray*}
This implies that the first assertion of Theorem \ref{mv}.
The variance of $\hat{\eta}_{D,\gamma}(x_D)=\vec{B}(x_D)^T\hat{\vec{b}}_D$ can be written as 
\begin{eqnarray*}
V[\hat{\eta}_{D,\gamma}(x_D)|\vec{X}_n]
&=&
\vec{B}(x_D)^T\Lambda_{D,\gamma}^{-1}V[G_D(\vec{b}_0,\lambda_n,\gamma_n)|\vec{X}_n]\Lambda_{D,\gamma}^{-1}\vec{B}(x_D)(1+o_P(1)).
\end{eqnarray*}
since it is easy to find that the conditional variance of $v_n(\vec{b}_0)$ can be shown to be $o_P(K_n/n)$. 
By noting
\begin{eqnarray*}
V[G_D(\vec{b}_0,\lambda_n,\gamma_n)|\vec{X}_n]
=
\frac{1}{n^2}Z_D^T V[\vec{y}|\vec{X}_n]Z_D=\frac{1}{n^2}Z_D^T W Z_D=\frac{1}{n}G_{D} +o_P((K_n/n)^{-1}\vec{1}\vec{1}^T),
\end{eqnarray*}
we have the second assertion as
\begin{eqnarray*}
&&V[\hat{\eta}_{D,\gamma}(x_D)|\vec{X}_n]\\
&&=
\frac{1}{n}\vec{B}(x_D)^T\Lambda_{D,\gamma}^{-1}G_{D} \Lambda_{D,\gamma}^{-1}\vec{B}(x_D)(1+o_P(1))\\
&&=
\frac{1}{n}\vec{B}(x_D)^T\Gamma_D(\lambda_{Dn})^{-1} \Gamma_D(0) \Gamma_D(\lambda_{Dn})^{-1}\vec{B}(x_D)(1+o_P(1)).
\end{eqnarray*}
Also it is easily confirmed by straightforward calculation with Lemma \ref{G1} that for $j\not= k$,
\begin{eqnarray*}
Cov(\hat{\eta}_j(x_j),\hat{\eta}_{k}(x_k))
=
\frac{1}{n}\vec{B}(x_j)^T\Lambda_{j,\gamma}^{-1}\left(\frac{1}{n}Z_j^T W Z_k\right)\Lambda_{k,\gamma}^{-1}\vec{B}(x_k)(1+o_P(1))
=O_P(n^{-1}),
\end{eqnarray*}
this completely the proof.
$\Box$

\vspace{5mm}

\subsection*{proof of Theorem \ref{cltpql}}

Let $\hat{\vec{b}}_P=[\hat{\vec{b}}_{1,P}^T\ \cdots\ \hat{\vec{b}}_{D,P}^T]^T$ and let $\tilde{\vec{b}}=(\tilde{\vec{b}}_1^T\ \cdots\ \tilde{\vec{b}}_D^T)^T$ be the maximizer of 
\begin{eqnarray*}
G(\vec{b},\Sigma_u,\tilde{\gamma}_n)=\frac{1}{n}\{\vec{y}^T (Z\vec{b})-\vec{1}^T c(Z\vec{b})\}+\frac{1}{n}\vec{1}^T h(\vec{y})-\frac{1}{2n} \vec{b}^T Q_{p+1}(\Sigma_u) \vec{b}-\frac{\tilde{\gamma}_n}{2n}\vec{b}^T \vec{b}
\end{eqnarray*}
with respect to $(\vec{b}_1^T\ \cdots\ \vec{b}_D^T)^T$, where $\tilde{\gamma}_n=o(K_n^{p-1}/\sigma_j^2)$ and let $\tilde{\eta}_j(x_j)=\vec{B}(x_j)^T\tilde{\vec{b}}_j (j=1,\cdots,D)$. 
Then the asymptotic normality of $[\tilde{\eta}_1(x_1)\ \cdots\ \tilde{\eta}_D(x_D)]^T$ can be obtained by the same manner to the proof of Theorem \ref{clt} with $m=p+1$. 
Similar to Lemma \ref{ex}, $\hat{\vec{b}}_P$ can be written as 
\begin{eqnarray*}
\hat{\vec{b}}_P-\tilde{\vec{b}}&=&(Z^T\tilde{W}Z+Q_{p+1}(\Sigma_u)+\tilde{\gamma}_nI)^{-1}G(\hat{\vec{b}}_P,\Sigma_u,\tilde{\gamma}_n)+r_n\\
&=&
\tilde{\gamma}_n(Z^T\hat{W}Z+Q_{p+1}(\Sigma_u)+\tilde{\gamma}_nI)^{-1}\hat{\vec{b}}_P+r_n
\end{eqnarray*}
where $\hat{W}=\diag[c^{\prime\prime}(Z\hat{\vec{b}}_P)]$ and $r_n=o_P(\sqrt{K_n/n}\vec{1})$ is the remainder.
Then Lemma \ref{G2} yields $\hat{\vec{b}}_P-\tilde{\vec{b}}=O_P(\tilde{\gamma}_nK_nn^{-1}\vec{1})=o_P(\sqrt{K_n/n}\vec{1})$, by which $\hat{\eta}_{j,P}(x_j)-\tilde{\eta}_j(x_j)=o_P(\sqrt{K_n/n}) (j=1,\cdots,D)$. 
This leads to Theorem \ref{cltpql}.
$\Box$

\def\bibname{References}

\end{document}